\documentclass[10pt,a4paper]{amsart}

\usepackage{graphicx,latexsym,amssymb}
\usepackage[all]{xy}
\CompileMatrices

\theoremstyle{plain}
\newtheorem{theorem*}{Theorem}
\newtheorem{theorem}{Theorem}[section]
\newtheorem{lemma}{Lemma}[section]
\newtheorem{proposition}{Proposition}[section]
\newtheorem{corollary}{Corollary}[section]
\newtheorem{corollary*}{Corollary}
\theoremstyle{definition}
\newtheorem{definition}{Definition}[section]
\theoremstyle{remark}
\newtheorem{remark}{Remark}[section]
\newtheorem{example}{Example}[section]
\newtheorem{claim}{Claim}[section]

\newtheorem{problem}{Problem}[section]
\numberwithin{equation}{section}
\numberwithin{figure}{section}

\font\nb=msbm10
\font\nbi=msbm8 at 7pt
\def \R{\hbox{\nb R}}

\def \Z{\hbox{\nb Z}}
\def \Zi{\hbox{\nbi Z}}
\def \Q{\hbox{\nb Q}}
\def \R{\hbox{\nb R}}
\def \Qi{\hbox{\nbi Q}}
\def \C{\hbox{\nb C}}

\def \Hom{\hbox{Hom}}
\def \Spin{\hbox{Spin}}
\def \Spinc{\hbox{Spin}^c}
\def \SO{\hbox{SO}}
\def \U{\hbox{U}}
\def \SU{\hbox{SU}}
\def \Spincscript{\hbox{\scriptsize Spin}^c}
\def \Spinscript{\hbox{\scriptsize Spin}}
\def \Uscript{\hbox{\scriptsize U}}
\def \SOscript{\hbox{\scriptsize SO}}
\def \Quad{\hbox{Quad}}

\def \Parall{\hbox{Parall}}

\def \spinc{{\rm{Spin}}$^c$}
\def \spin{{\rm{Spin}}}

\begin{document}

\title[Quadratic functions and complex spin structures on  $3$--manifolds]
{Quadratic functions and complex spin structures on  three--manifolds}
\author[F. Deloup]{Florian Deloup}
\author[G. Massuyeau]{Gw\'ena\"el Massuyeau}
\begin{abstract}
We show how the space of complex spin structures of a closed oriented three--manifold
embeds naturally into a space of quadratic functions associated to its linking pairing.
Besides, we extend the Goussarov--Habiro theory of finite type invariants
to the realm of compact oriented three--manifolds equipped with a complex spin structure.
Our main result states that two closed oriented three--manifolds endowed with a complex spin structure are
undistinguishable by complex spin invariants of degree zero if, and only if, their associated quadratic functions are isomorphic.
\end{abstract}
\maketitle

Complex spin structures, or \spinc--structures, 
are additional structures with which manifolds may be  equipped. 
They are needed to define the Seiberg--Witten invariants of $4$--manifolds, 
as well as the Heegaard--Floer homologies of $3$--manifolds by Ozsv\'ath and Szab\'o.
Any closed oriented $3$--manifold $M$ can be endowed with a \spinc--structure and, in that case,
\spinc--structures are in canonical correspondence with Euler structures.
The latter are classes of nonsingular vector fields on $M$ 
which have been introduced by Turaev in order to refine Reidemeister torsion.

In this paper, we investigate the r\^ole played by quadratic functions in the topology of closed
oriented $3$--manifolds equipped with a \spinc--structure or, equivalently, an Euler structure.\\

Extending constructions from \cite{LL,MS,LW},
we associate, to any closed oriented $3$--manifold $M$ with a \spinc--structure
$\sigma$, its \emph{linking quadratic function}
$$
\xymatrix{
{H_{2}(M;\Q/\Z)} \ar[r]^-{\phi_{M,\sigma}} \ar[r] & {\Q/\Z}.
}
$$
The function $\phi_{M,\sigma}$ is quadratic in the sense that the symmetric pairing
defined by $(x,y) \mapsto \phi_{M,\sigma}(x+y) - \phi_{M,\sigma}(x) - \phi_{M,\sigma}(y)$
is bilinear. Moreover, this symmetric bilinear pairing 
coincides with $L_M:=\lambda_M\circ (B\times B)$ where 
$$
\xymatrix{ \hbox{Tors }H_{1}(M;\Z) \times \hbox{Tors }H_{1}(M;\Z) \ar[r]^-{ \lambda_{M}}& \Q/\Z}
$$
is the linking pairing of $M$ and $B$ denotes the Bockstein homomorphism associated to the short exact
sequence of coefficients $0\to \Z \to \Q \to \Q/\Z \to 0$. In contrast with $\phi_{M,\sigma}$, 
the bilinear pairing $L_M$ does not depend on $\sigma$.
\spinc--structures on a given manifold $M$ are determined by their corresponding quadratic functions.

\begin{theorem*}
\label{th:embedding} 
Let $M$ be a closed connected oriented $3$--manifold. 
The map $\sigma \mapsto \phi_{M,\sigma}$ 
defines a canonical embedding 
$$
\xymatrix{
{\rm{Spin}}^c(M)\ \ar@{^{(}->}[r]^{\phi_M} &{\rm{Quad}}(L_M)
}
$$
from the set of \spinc--structures on $M$
to the set of quadratic functions with $L_M$ as associated bilinear pairing.
\end{theorem*}
\noindent
Via the map $\phi_M$, topological notions can be put 
in correspondence with algebraic ones. For instance, 
the Chern class $c(\sigma)\in H^2(M)$ of the \spinc--structure $\sigma$
corresponds to the homogeneity defect $d_{\phi_{M,\sigma}}: H_2(M;\Q/\Z) \to \Q/\Z$
of the quadratic function $\phi_{M,\sigma}$, which is defined by
$d_{\phi_{M,\sigma}}(x)$ $=$ $\phi_{M,\sigma}(x)-\phi_{M,\sigma}(-x)$.

When the Chern class $c(\sigma)$ is torsion, $\phi_{M,\sigma}$
happens to factor through $B$ to a quadratic function
$$
\xymatrix{
{\hbox{Tors }H_{1}(M;\Z)} \ar[r]^-{\phi_{M,\sigma}} \ar[r] & {\Q/\Z}
}
$$ 
with $\lambda_M$ as associated bilinear pairing and is equivalent
to the quadratic function constructed by Looijenga and Wahl \cite{LW} (see also \cite{Gille,De}).
In particular, the \spinc--structure may arise from a classical spin structure, 
or Spin--structure. In that case, which is detected by the 
vanishing of $c(\sigma)$, the quadratic function $\phi_{M,\sigma}$ is homogeneous and coincides
with yet earlier constructions due to Lannes and Latour \cite{LL}, as well as 
Morgan and Sullivan \cite{MS} (see also \cite{TCohomology,KT}).\\

The linking quadratic function is used here to solve a problem related to the 
theory of finite type invariants by Goussarov and Habiro.
Their theory \cite{Goussarov,Habiro,GGP} deals with compact oriented $3$--manifolds
and is based on an elementary move called $Y$--surgery.
The $Y$--equivalence, which is defined to be the equivalence relation among
such manifolds generated by this move, has been characterized by Matveev in the closed case \cite{Matveev}.
This characterization amounts to recognize the degree $0$ invariants of the theory. 
His result, anterior to the work of Goussarov and Habiro, can be re-stated as follows:
two closed oriented $3$--manifolds $M$ and $M'$ are $Y$--equivalent if and only if they have 
isomorphic pairs (homology, linking pairing). A \spin--refinement of the Goussarov--Habiro theory 
(the possibility of which was announced in \cite{Goussarov} and \cite{Habiro}) 
has also been considered in \cite{Mas}, where Matveev's theorem is extended to 
closed oriented $3$--manifolds equipped with a Spin--structure.

We show that the $Y$--surgery move makes sense for closed oriented 
$3$--manifolds equipped with a \spinc--structure as well.
The equivalence relation generated by this move among such manifolds is called, here, 
\emph{$Y^c$--equivalence}. It follows that there exists a \spinc--refinement of the Goussarov--Habiro theory.
Our main result is a characterization of the $Y^c$--equivalence relation
in terms of the linking quadratic function.
In order to state this more precisely, let us fix a few notations.

Given an isomorphism $\psi: H_{1}(M;\Z) \to H_{1}(M';\Z)$, the dual
isomorphism to $\psi$ by the intersection pairings is denoted by 
$\psi^\sharp: H_2\left(M';\Q/\Z\right) \rightarrow H_2(M;\Q/\Z)$:
$$
\forall x\in H_1(M;\Z),\ \forall y'\in H_2(M';\Q/\Z),\quad
x\bullet \psi^\sharp(y')=\psi(x)\bullet y' \in \Q/\Z.
$$ 
Also, given sections $s$ and $s'$ of the surjections
$B:H_2\left(M;\Q/\Z\right)\to {\hbox{Tors }}H_1(M;\Z)$ and
$B:H_2\left(M';\Q/\Z\right)\to {\hbox{Tors }}H_1(M';\Z)$ respectively, 
we say that $s$ and $s'$ are $\psi$--compatible if the diagram 
$$\xymatrix{
H_2\left(M';\Q/\Z\right) \ar[d]_-{\psi^{\sharp}}^{\simeq} 
& \hbox{Tors }H_1(M';\Z) \ar[l]_-{s'} \\
H_2\left(M;\Q/\Z\right) 
& \hbox{Tors }H_1(M;\Z) \ar[l]^-{s} \ar[u]^-{\psi|}_{\simeq}
}$$
commutes. We denote by $P$ a Poincar\'e isomorphism and
we recall that the Gauss sum of a quadratic function $q:G\to \Q/\Z$, defined on a finite
Abelian group $G$, is the complex number $\sum_{x \in G} \exp(2 i \pi q(x))$. 
 
\begin{theorem*} 
\label{th:Matveev_spinc}
Let $(M,\sigma)$ and $(M',\sigma')$ be two closed 
connected oriented $3$--manifolds with \spinc--structure. 
The following assertions are equivalent:
\begin{enumerate}
\item[$(1)$] The \spinc--manifolds $(M,\sigma)$ and $(M',\sigma')$ 
are $Y^c$--equivalent.
\item[$(2)$] There is an isomorphism 
$\psi: H_{1}(M;\Z) \to H_{1}(M';\Z)$ 
such that $\phi_{M',\sigma'} = \phi_{M,\sigma} \circ \psi^{\sharp}$.
\item[$(3)$] There is an isomorphism 
$\psi: H_{1}(M;\Z) \to H_{1}(M';\Z) $ such that
\begin{itemize}
\item[$\centerdot$] $\lambda_{M} = \lambda_{M'} \circ \left(\psi| \times \psi| \right)$,
\item[$\centerdot$] $\psi\left(P^{-1}c(\sigma)\right) = P^{-1}c(\sigma')$, 
\item[$\centerdot$] for some $\psi$--compatible sections $s$ and $s'$ 
of the Bockstein homomorphisms, $\phi_{M,\sigma}\circ s$ and $\phi_{M',\sigma'}\circ s'$ 
have identical Gauss sums. 
\end{itemize}
\end{enumerate}
\end{theorem*}
 
Two special cases deserve to be singled out. 
First, consider manifolds whose first homology group is torsion free.
The following result is deduced from Theorem \ref{th:Matveev_spinc}.
\begin{corollary*}
\label{cor:torsion_free}
Let $(M,\sigma)$ and $(M',\sigma')$ be two closed connected oriented 
$3$--manifolds with \spinc--structure, such that $H_1(M;\Z)$ and $H_1(M';\Z)$
are torsion free. The following assertions are equivalent:
\begin{enumerate}
\item[$(1)$] The \spinc--manifolds $(M,\sigma)$ and $(M',\sigma')$ 
are $Y^c$--equivalent.
\item[$(2)$] There is an isomorphism 
$\psi: H_{1}(M;\Z) \to H_{1}(M';\Z)$ such that
$\psi\left(P^{-1}c(\sigma)\right)$$ = P^{-1}c(\sigma')$. 
\end{enumerate}
\end{corollary*}
\noindent
Second, consider the case of rational homology $3$--spheres. 
According to what has been said above, 
if $M$ is an oriented rational homology $3$--sphere, 
then $\phi_{M,\sigma}$ can be regarded 
as a quadratic function $H_{1}(M;\Z)\rightarrow \Q/\Z$ with $\lambda_{M}$ as associated bilinear pairing. 
In that case, Theorem \ref{th:Matveev_spinc} specializes to the following corollary.
\begin{corollary*}
\label{cor:QHS}
Let $(M,\sigma)$ and $(M',\sigma')$ be two oriented rational homology 
$3$--spheres with \spinc--structure. The following assertions are equivalent:
\begin{enumerate}
\item[$(1)$] The \spinc--manifolds $(M,\sigma)$ and $(M',\sigma')$ 
are $Y^c$--equivalent.
\item[$(2)$] There is an isomorphism 
$\psi: H_{1}(M;\Z) \to H_{1}(M';\Z)$ such that $\phi_{M,\sigma} = \phi_{M',\sigma'} \circ \psi$.
\item[$(3)$] There is an isomorphism $\psi: H_{1}(M;\Z) \to H_{1}(M';\Z) $ such that
\begin{itemize}
\item[$\centerdot$] $\lambda_{M} = \lambda_{M'} \circ \left(\psi\times \psi\right)$,
\item[$\centerdot$] $\psi\left(P^{-1}c(\sigma)\right) = P^{-1}c(\sigma')$, 
\item[$\centerdot$] $\phi_{M,\sigma}$ and $\phi_{M',\sigma'}$ have identical Gauss sums.
\end{itemize}
\end{enumerate}
\end{corollary*}
\vspace{0.5cm}

The paper is organized as follows. In Section \ref{sec:spinc}, we briefly review
\spinc--structures from a general viewpoint. Next, we restrict ourselves to the dimension $3$,
in which case one can work with Euler structures as well.
At the end of the section, the technical problem of gluing Spin$^c$--structures is considered. 
This is needed to define the $Y$--surgery move in the setting of manifolds equipped 
with a \spinc--structure, since this move is defined as a ``cut and paste" operation.
Our gluing lemma involves \spinc--structures, on a compact oriented $3$--manifold with boundary,
which are relative to a fixed \spin--structure on the boundary. 

Section \ref{sec:quad_spinc} is devoted to the construction and study
of the linking quadratic function. First, we give a combinatorial description 
of the \spinc--structures of a given closed oriented $3$--manifold
presented by surgery along a link in $\mathbf{S}^3$. 
This leads to a \spinc--refinement of Kirby's theorem.
Next, we define the quadratic function $\phi_{M,\sigma}$ associated
to a closed $3$--dimensional \spinc--manifold $(M,\sigma)$: 
this is done essentially by defining a cobordism invariant 
of singular $3$--dimensional \spinc--manifolds over $K(\Q/\Z,1)$.
The quadratic function $\phi_{M,\sigma}$ can be computed combinatorially
as soon as $(M,\sigma)$ is presented by surgery along a link in $\mathbf{S}^3$.
We prove Theorem \ref{th:embedding} and some other basic properties of the map $\phi_M$. 
Lastly, regarding $\sigma$ as an Euler structure,  we give for $\phi_{M,\sigma}$ 
an intrinsic formula that does not make reference to the dimension $4$ anymore. 
This is obtained by presenting, \`a la Sullivan, elements of $H_2(M;\Q/\Z)$ 
as immersed surfaces with $n$--fold boundary.

In Section \ref{sec:FTI}, the $Y^c$--surgery move is defined
using the above mentioned gluing lemma.
Next, Theorem \ref{th:Matveev_spinc} is proved working with surgery presentations of \spinc--manifolds. 
We use the material of the previous section and a result due to Matveev, Murakami and Nakanishi 
\cite{Matveev,MN} on ordered oriented framed links having the same linking matrix. 
Some algebraic ingredients about quadratic functions on torsion Abelian groups are needed as well. 
Those results, some of them well--known in the case of finite Abelian groups, have been proved aside in \cite{DM1}. 
We conclude this paper by giving some applications of Theorem \ref{th:Matveev_spinc} 
and stating some problems.\\

\noindent
\textbf{Acknowledgments.} 
F. D. has been supported by a Marie Curie Research Fellowship (HPMF-CT-2001-01174).
G. M. thanks G\'erald Gaudens for useful discussions, and
Christian Blanchet for his help and encouragements.

\section{Complex spin structures on three--manifolds}

\label{sec:spinc}
In this section, we review \spinc--structures and other related structures,
with special emphasis on the dimension $3$.
We also give a gluing lemma for \spinc--structures.

\subsection{Some conventions}

\label{subsec:miscellaneous}
In this paper, any manifold $M$ is assumed to be compact, smooth and oriented.
We denote by $-M$ the manifold obtained from $M$ by reversing its orientation.
If $M$ has non-empty boundary, $\partial M$ has the orientation given 
by the ``outward normal vector first'' rule. 
The oriented tangent bundle of $M$ is denoted by T$M$.

Vector bundles will be stabilized from the left side. A section of a vector bundle
is said to be \emph{nonsingular} if it does not vanish at any point.

If $G$ is an Abelian group, a \emph{$G$--affine space} $A$ is a set $A$
on which $G$ acts freely and transitively. The affine action is denoted additively; thus,
for $a$, $a'\in A$, the unique element $g\in G$ satisfying $a'=a + g$ will be written $a' - a$.

Unless otherwise specified, all (co)homology groups are assumed 
to be computed with integer coefficients.

\subsection{Complex spin structures}

In this subsection, we consider a $n$--manifold $M$.
We recall basic facts about Spin$^c$--structures on $M$,
adopting a viewpoint which is analogous to that used in \cite{BM} for \spin--structures.

\subsubsection{From  ${\rm{Spin}}^c$ onto ${\rm{SO}}$}

\label{subsubsec:up}
Let $n\geq 1$ be an integer. The group Spin$(n)$ is the $2$--fold covering
of the special orthogonal group $\SO(n)$: 
$$
\xymatrix{
{1}  \ar[r] & {\Z_2} \ar[r] & {\Spin(n)} \ar[r] & \SO(n) \ar[r] & 1.
}
$$
The group ${\Spin}^c(n)$ is defined by
$$
{\Spin}^c(n)=\frac{{\Spin(n)}\times \U(1)}{\Z_2}
$$
where $\Z_2$ is generated by $\left[(-1,-1)\right]$, 
hence the following short exact sequence of groups:
$$\xymatrix{
1 \ar[r] & \U(1) \ar[r] & {\Spin}^c(n) \ar[r]^-{\pi} & \SO(n) \ar[r] & 1
}
$$
where the first map sends $z$ to $\left[(1,z)\right]$
and where $\pi$ is induced by the projection of Spin$(n)$ onto SO$(n)$.

The inclusion of $\SO(n)$ into $\SO(n+1)$, defined by $A\mapsto (1) \oplus A$, induces a monomorphism
$\xymatrix{\Spinc(n)\ \ar@{>->}[r] & \Spinc(n+1)}$ such that the diagram 
$$
\xymatrix{
\Spinc(n)\ \ar@{>->}[r] \ar@{->>}[d]_-{\pi} & \Spinc(n+1) \ar@{->>}[d]^-{\pi}\\
\SO(n)\ \ar@{>->}[r] & \SO(n+1).
}
$$
commutes, hence a diagram at the level of classifying spaces:
\begin{equation}
\label{eq:strict}
\xymatrix{
\hbox{B}\Spinc(n) \ar[r] \ar[d]_-{\hbox{B} \pi} & \hbox{B} \Spinc(n+1) \ar[d]^-{\hbox{B}\pi}\\
\hbox{B}\SO(n) \ar[r] & \hbox{B}\SO(n+1).
}
\end{equation}
Here, we take B$\SO(n)$ to be the Grassman manifold of oriented $n$--planes in $\R^\infty$ 
and the map $\hbox{B}\SO (n) \to \hbox{B}\SO (n+1)$ to be the usual one.
We fix the classifying spaces $\hbox{B}\Spinc(n)$ (in their homotopy equivalence classes)
and, next,  we fix the maps $\hbox{B} \pi :\hbox{B}\Spinc(n) \to \hbox{B}\SO (n)$ (in their homotopy classes) 
to be fibrations. Then, the map from $\hbox{B}\Spinc(n)$ to $\hbox{B}\Spinc(n+1)$ 
is choosen (in its homotopy class) to make diagram (\ref{eq:strict}) \emph{strictly} commute.

We denote by $\gamma_{{\SOscript}(n)}$ the universal $n$--dimensional oriented vector bundle
over B$\SO(n)$. Let $\gamma_{{\Spinscript}^c(n)}$ be the pull--back of $\gamma_{{\SOscript}(n)}$ by B$\pi$.
Thanks to (\ref{eq:strict}), there is a well--defined morphism between $(n+1)$--dimensional oriented vector bundles
$\R \oplus \gamma_{\Spincscript(n)} \to \gamma_{\Spincscript(n+1)}$ induced by the usual one
$\R \oplus \gamma_{\SOscript(n)} \to \gamma_{\SOscript(n+1)}$. 

\subsubsection{Rigid ${\rm{Spin}}^c$--structures}

\label{subsubsec:def_spinc}

Recall that $M$ is a $n$--manifold to which some conventions, stated
in \S \ref{subsec:miscellaneous}, apply.
\begin{definition}
\label{def:def_spinc}
A \emph{rigid \spinc--structure} on $M$ is a morphism 
$\hbox{T}M \to \gamma_{{\Spinscript}^c(n)}$ between $n$--dimensional oriented vector bundles. 
A \emph{\spinc--structure} (or \emph{complex spin structure}) 
on $M$ is a homotopy class of rigid 
\spinc--structures on $M$. We denote by ${\Spin}^c_r(M)$ 
the set of rigid \spinc--structures on $M$, and by ${\Spin}^c(M)$
the set of its \spinc--structures.
\end{definition}

\noindent
Obviously, a different choice of the classifying space BSpin$^c(n)$ (in its homotopy type)  
or a different choice of the map B$\pi$ (in its homotopy class)  would
lead to a different notion of rigid \spinc--structure, but would
not affect the definition of a \spinc--structure. 
Rigid structures will be used later to define gluing maps.

Let $\beta$ be the Bockstein homomorphism
associated to the short exact sequence of coefficients
\begin{displaymath}
\xymatrix{
0 \ar[r] & {\Z} \ar[r]^{\cdot 2} & {\Z}
\ar[r] & {\Z_2} \ar[r] & 0.
}
\end{displaymath}
The fibration $\hbox{B}\pi:\hbox{B} \hbox{Spin}^c(n) \to \hbox{B}\SO(n)$ has fiber $\hbox{B}\U(1)\simeq K(\Z,2)$
and, indeed, is a principal fibration with 
characteristic class $w:=\beta w_2 \in H^3(\hbox{B}\SO(n))$, where $w_2$ is the second
Stiefel--Whitney class. Then, by obstruction theory, we obtain the following 
well--known fact about existence and parametrization of Spin$^c$--structures.
\begin{proposition}
\label{prop:absolute_obstruction}
The manifold $M$ can be given a \spinc--structure if and only if the 
cohomology class $\beta w_2(M) \in H^3(M)$ vanishes.
In that case, ${\rm{Spin}}^c(M)$ is a $H^2(M)$--affine space.
\end{proposition}

One may easily verify that the homotopy--theoretical definition of a \spinc--structure,
which we have adopted here, agrees with the usual one.
\begin{lemma} 
\label{lem:usual_definition}
Suppose that $M$ is equipped with a Riemannian metric and
denote by ${\rm{SO}}({\rm{T}}M)$ the bundle of its oriented orthonormal frames.
A \spinc--structure on $M$ is equivalent to an isomorphism class of pairs $(\eta,H)$, 
where $\eta$ is a principal ${\rm{Spin}}^c(n)$--bundle over $M$
and where $H:\eta/{\rm{U}}(1) \to {\rm{SO}}({\rm{T}}M)$ is a principal {\rm{SO}}$(n)$--bundle isomorphism.
\end{lemma}

To go to the point, we have only defined (rigid) \spinc--structures on the manifold $M$.
Nevertheless,  the notion of a (rigid) \spinc--structure obviously 
extends to any oriented vector bundle over any base space.

\begin{remark}
\label{rem:stable}
Thanks to the map $\R \oplus \gamma_{\Spincscript(n)} \to \gamma_{\Spincscript(n+1)}$
constructed at the end of \S \ref{subsubsec:up}, a rigid \spinc--structure
on $\hbox{T} M$ gives rise to one on $\R \oplus \hbox{T} M$. This induces a canonical map
$$
 \hbox{Spin}^c(M)=\hbox{Spin}^c(\hbox{T} M) \to \hbox{Spin}^c(\R \oplus \hbox{T} M)
$$
which is $H^2(M)$--equivariant and, so, bijective. 
Thus, a \spinc--structure on $M$ is equivalent to a \spinc--structure on its stable oriented tangent bundle.
\end{remark}

\subsubsection{Orientation reversal}

The \emph{time--reversing} map is the orientation--reversing automorphism 
of $\R\oplus \hbox{T} M$ defined by $(t,v)\mapsto (-t,v)$. Composition with that map
transforms a rigid \spinc--structure on $\R\oplus \hbox{T} M$ to one on $\R\oplus \hbox{T} (-M)$.
So, by Remark \ref{rem:stable}, we get a canonical $H^2(M)$--equivariant map 
$$ 
\xymatrix{
{\hbox{Spin}^c(M)} \ar[r]^-{-} & {\hbox{Spin}^c(-M)}.
}
$$

\subsubsection{Relative \spinc--structures}

\label{subsubsec:relative_Spinc}
Suppose that $M$ has some boundary and  fix a rigid 
structure $s\in \hbox{Spin}^c_r\left(\hbox{T}M|_{\partial M}\right)$ over $\partial M$.
\begin{definition}
A \emph{\spinc--structure on $M$ relative to $s$} 
is a homotopy class rel $\partial M$ of rigid \spinc--structures on $M$ that extend $s$. We denote
by $\hbox{Spin}^c(M,s)$ the set of such structures.
\end{definition}
\noindent

The following relative version of Proposition \ref{prop:absolute_obstruction}
is also proved by obstruction theory applied to the fibration B$\pi$.
\begin{proposition}
\label{prop:relative_obstruction}
There exists a rigid \spinc--structure on $M$ that extends $s$ if and only if a certain cohomology class
\begin{displaymath} 
w(M,s)\in H^3(M,\partial M) 
\end{displaymath}
vanishes. In that case, ${\rm{Spin}}^c(M,s)$ is a $H^2(M,\partial M)$--affine space.
\end{proposition}

\subsubsection{Restriction to the boundary}

\label{subsubsec:boundary}
Suppose that $M$ has some boundary.
Observe that there is a well--defined homotopy class of isomorphisms between
the oriented vector bundles $\R \oplus \hbox{T} \partial M$ and $\hbox{T}M|_{\partial M}$, which is
defined by any section of T$M|_{\partial M}$ transverse to $\partial M$ and directed outwards.

In particular, a \spinc--structure on $\hbox{T}M|_{\partial M}$ can be identified without ambiguity
to a \spinc--structure on $\partial M$.
Thus, we get a canonical \emph{restriction} map
$$
\xymatrix{
\hbox{Spin}^c(M) \ar[r] & \hbox{Spin}^c(\partial M)
}
$$
which is affine over the homomorphism $H^2(M) \to H^2(\partial M)$ induced by inclusion.

\subsubsection{From ${\rm{Spin}}$  to ${\rm{Spin}}^c$} 

\label{subsubsec:pont}

Proceeding as in \S \ref{subsubsec:def_spinc}, we define the set Spin$_r(M)$ 
of \emph{rigid \spin--structures} on $M$ and the set Spin$(M)$ of \emph{\spin--structures} on $M$.
The latter is a $H^2(M;\Z_2)$--affine space as soon as $w_2(M)$ vanishes. 
The reader is refered to \cite{BM} for details.\footnote{In \cite{BM}, rigid \spin--structures are called 
``$w_2$--structures'' and are defined on the stable oriented tangent bundle. 
An observation similar to that given in Remark \ref{rem:stable} for \spinc--structures applies
to \spin--structures.} The group homomorphism
$$
\xymatrix{{\hbox{Spin}}(n) \ar[r]^{\beta} &{\hbox{Spin}^c}(n)}
$$
defined by $\beta(x)=\left[(x,1)\right]$, makes the two projections onto SO$(n)$ agree.
This allows us to define a morphism $\gamma_{\Spinscript(n)} \to \gamma_{\Spinscript^c(n)}$ between oriented 
$n$--dimensional vector bundles, the composition with which transforms a rigid \spin--structure $u$
to a rigid \spinc--structure denoted by $\beta(u)$. Thus, we get a canonical map
$$
\xymatrix{
\hbox{Spin}(M) \ar[r]^{\beta} & \hbox{Spin}^c(M)
}
$$
which is affine over the Bockstein homomorphism 
$\beta:H^1(M;\Z_2) \to H^2(M)$. 

If $M$ has some boundary, we define \emph{relative \spin--structures} on $M$ as well.
Their construction goes as in \S \ref{subsubsec:relative_Spinc}. 
Thus, for a fixed $s\in \hbox{Spin}_r\left(\hbox{T}M|_{\partial M}\right)$, we get a map
$$
\xymatrix{
\hbox{Spin}(M, s) \ar[r]^-{\beta} & \hbox{Spin}^c(M, \beta s)
}
$$
which is affine over the Bockstein homomorphism $\beta:H^1(M,\partial M;\Z_2) \to H^2(M, \partial M)$. 

\subsubsection{From ${\rm{U}}$  to ${\rm{Spin}}^c$} 

\label{subsubsec:U}
Let $m$ be an integer such that $n\leq 2m$. 
We take BU$(m)$ to be the Grassman manifold 
of complex $m$--planes in $\C^\infty$. The map BU$(m) \to \hbox{B}\SO(2m)$, 
which consists in forgetting the complex structure on a complex $m$--plane,
represents the usual inclusion of $\U(m)$ into $\SO(2m)$.
We define $\gamma_{\Uscript(m)}$ to be the pull--back of $\gamma_{\SOscript(2m)}$ 
by this map BU$(m) \to \hbox{B}\SO(2m)$, 
which can be identified with the $2m$--dimensional oriented vector bundle underlying 
the universal $m$--dimensional complex vector bundle. Then, as we did in the Spin 
and Spin$^c$ cases, we could define a ``rigid U--structure'' on $ \R^{2m-n} \oplus \hbox{T} M$ to be 
a morphism $\R^{2m-n} \oplus \hbox{T}M \to \gamma_{\Uscript(m)}$ 
between $2m$--dimensional oriented vector bundles.
Such a morphism induces a complex structure on $\R^{2m-n} \oplus \hbox{T}M$
by pulling back the canonical one on $\gamma_{\Uscript(m)}$
and, conversely, any complex structure on $\R^{2m-n} \oplus \hbox{T}M$
inducing the given orientation arises that way.
Then, a ``U--structure'' on  $\R^{2m-n} \oplus \hbox{T}M$ is equivalent to 
a homotopy class of complex structures on $\R^{2m-n} \oplus \hbox{T}M$ compatible with the given orientation. 

There is a canonical way to embed U$(m)$ into Spin$^c(2m)$: see, for instance, \cite[Proposition D.50]{GGK}.
This inclusion
$$
\xymatrix{{\hbox{U}}(m)\ \ar@{>->}[r]^-{\omega} & {\hbox{Spin}^c}(2m)}
$$
makes the two maps to SO$(2m)$ commute. This allows us to define a morphism 
$\gamma_{\Uscript(m)} \to \gamma_{\Spinscript^c(2m)}$ between oriented $2m$--dimensional vector bundles, 
the composition with which transforms a ``rigid U--structure'' on $\R^{2m-n} \oplus \hbox{T}M$
to a rigid \spinc--structure on it. As a consequence of Remark \ref{rem:stable}, we get a canonical map
$$
\xymatrix{\hbox{U}^{s}(M) \ar[r]^-\omega &  {\hbox{Spin}^c}(M)}
$$
from the set of stable complex structures on T$M$ compatible with the orientation 
to the set of Spin$^c$--structures on $M$.
(See \cite[Proposition D.57]{GGK} for a construction of $\omega$ 
involving the usual definition of a \spinc--structure.) 

\subsubsection{Chern class}

\label{subsubsec:Chern}

A \spinc--structure $\alpha$ on  $M$ induces an isomorphism class of principal Spin$^c$(n)--bundles
over $M$ and, so, an  isomorphism class of principal U$(1)$--bundles 
thanks to the homomorphism $\hbox{Spin}^c(n)\to \U(1)$ defined by $[(x,y)]\mapsto y^2$. 
The first Chern class of the latter is denoted by $c(\alpha)$. We get a \emph{Chern class} map
$$
\xymatrix{\Spinc(M) \ar[r]^-c & H^2(M)}
$$
which is affine over the doubling map defined by $x\mapsto 2x$.
When $c(\alpha)$ belongs to Tors $H^2(M)$, the \spinc--structure $\alpha$ is said to be \emph{torsion}.

\subsection{Complex spin structures in dimension $3$}

In this subsection, we turn to $3$--manifolds which, by \S \ref{subsec:miscellaneous}, 
are assumed to be compact smooth and oriented.
The preliminary remark is that any $3$--manifold $M$ can be endowed 
with a \spinc--structure, since $w_2(M)$ is well--known to vanish.

We start by removing the rigidity of relative \spinc--structures
which is still remaining along the boundary. Next, we recall Turaev's observation
that \spinc--structures can be regarded as classes of vector fields. 
This holds true in the relative case as well.

\subsubsection{Relative \spinc--structures}

Let $M$ be a $3$--manifold with boundary and let $\sigma$ be a \spin--structure on $\partial M$.
We define  \spinc--structures on $M$ which are relative to $\sigma$. 
Note that, thanks to the observation initiating \S \ref{subsubsec:boundary},
one can identify $\sigma\in \Spin(\partial M)$ to a \spin--structure on $\hbox{T}M|_{\partial M}$.
\begin{lemma}
\label{lem:relative_spinc}
For any rigid \spin--structure $s$ on ${\rm{T}}M|_{\partial M}$ representing $\sigma$ 
(which we denote by $s\in \sigma$), the rigid \spinc--structure $\beta(s)$ can be extended to $M$.
Moreover, for any $s,s' \in \sigma$, there exists a canonical $H^2(M,\partial M)$--equivariant bijection
$$
\xymatrix{
{{\rm{Spin}}^c\left(M,\beta s\right)} \ar[r]^-{\rho_{s,s'}} & {{\rm{Spin}}^c(M,\beta s')}.
}
$$
Lastly, for any $s,s',s'' \in \sigma$, we have that $\rho_{s',s''} \circ \rho_{s,s'} = \rho_{s,s''}$.
\end{lemma}
\begin{definition}
\label{def:relative_spinc}
A \emph{\spinc--structure on $M$ relative to $\sigma$} is a pair $(u,s)$ where $s\in \sigma$ 
and $u\in {{\rm{Spin}}^c\left(M,\beta s\right)}$,
two such pairs $(u,s)$ and $(u',s')$ being considered as equivalent when $u'=\rho_{s,s'}(u)$.
The set of such structures is denoted by ${\Spin}^c(M,\sigma)$ and can naturally be
given the structure of a $H^2(M,\partial M)$--affine space. 
\end{definition}

\begin{remark}
\label{rem:relative}
There is an analogue to Lemma \ref{lem:usual_definition} that formulates what
a \spinc--structure on $M$ relative to $\sigma$  is in terms of principal bundles.
\end{remark}

\begin{example} 
\label{ex:tores_spinc}
Suppose that $\partial M$ is a disjoint union of tori. 
The $2$--torus has a distinguished Spin--structure $\sigma^0$ that is
induced by its Lie group structure. Using the previous remark,
it can be verified that a \spinc--structure on $M$ relative to the union of copies of $\sigma^{0}$ 
is equivalent to a relative \spinc--structure in the sense of Turaev \cite[\S 1.2]{TSW}.
\end{example}

\begin{proof}[Proof of Lemma \ref{lem:relative_spinc}]
Let $w_2(M,s)\in H^2(M,\partial M;\Z_2)$ denote the obstruction
to extend $s$ to a rigid \spin--structure on $M$. We have that
$$
\beta(w_2(M,s))=w(M,\beta s) \in H^3(M,\partial M).
$$
Thus, $w(M,\beta s)$ is of order at most $2$ and, so, vanishes.

We now prove the second statement. Let $\varphi:[-1,0] \times \partial M \hookrightarrow M$
be a collar neighborhood of $\partial M$. In particular,
$\varphi$ induces a specific isomorphism between $\R \oplus \hbox{T} \partial M$
and $\hbox{T}M|_{\partial M}$: the rigid \spin--structures on $\R \oplus \hbox{T} \partial M$
corresponding to $s$ and $s'$ are denoted by $s_0$ and $s_1$ respectively.
By assumption, $s_0$ and $s_1$ are homotopic: let $S=(s_t)_{t\in [0,1]}$ be
such a homotopy. This defines a rigid \spin--structure $S$ on
$[0,1] \times \partial M$ by identifying, at each time $t$,
$\R \oplus \hbox{T} \partial M$ with the restriction of
$\hbox{T}([0,1] \times \partial M)$ to $t \times \partial M$.
The same collar neighborhood allows us to 
define a smooth gluing $M \cup \left([0,1] \times \partial M\right)$,
as well as a positive diffeomorphism $\tilde{\varphi}: M \to M \cup \left([0,1] \times \partial M\right)$
(based on the affine identification between $[-1,0]$ and $[-1,1]$). Consider the map
$$
\xymatrix{ {\rm{Spin}}^c(M,\beta s) \ar[r]^-{\rho_{S}}& {\rm{Spin}}^c(M,\beta s')}
$$
defined, for any $u \in {\rm{Spin}}^c_r(M)$ extending $\beta(s)$, by 
$\rho_S([u])=\left[\left(u \cup \beta(S) \right) \circ \hbox{T} \tilde{\varphi}\right]$.

The map $\rho_{S}$ is $H^2(M,\partial M)$--equivariant and is independent of the choice of $\varphi$. 
So, we are left to prove that $\rho_S$ does not depend on the choice
of the homotopy $S$ between $s_0$ and $s_1$, which will allow us to set $\rho_{s,s'}= \rho_S$. 
To see that, consider the map $\beta$ constructed in \S \ref{subsubsec:pont}
from $\hbox{Spin}\left([0,1] \times \partial M,0 \times (- s_0) \cup 1 \times s_1 \right)$
to $\hbox{Spin}^c\left([0,1] \times \partial M,0 \times (-\beta s_0) \cup 1 \times (\beta s_1)\right)$,
where $-s_0 \in \hbox{Spin}_r(\R \oplus \hbox{T} (-\partial M))$ is obtained from $s_0$ by time--reversing.
The Bockstein homomorphism $\beta$ from $H^1([0,1] \times \partial M, \partial [0,1] \times \partial M;\Z_2)$
to $H^2([0,1] \times\partial M,\partial [0,1] \times\partial M)$  is trivial, 
since its codomain is isomorphic to the free Abelian group $H_1(\partial M)$. 
It follows that the former map $\beta$ collapses, and the conclusion follows.
\end{proof}
\begin{remark}
\label{rem:from_spin_to_spinc_relative}
The set of \emph{\spin--structures on $M$ relative to $\sigma$} is defined to be
\begin{displaymath}
\Spin(M,\sigma)=\left\{\alpha \in \Spin(M)\ 
:\ \alpha|_{\partial M}=\sigma \right\},
\end{displaymath}
which may be empty. One can construct a canonical map
$$
\xymatrix{
{\Spin(M,\sigma)}\ar[r]^\beta & {\Spin^c(M,\sigma)}
}
$$
by means of a rigid \spin--structure $s$ on ${\rm{T}}M|_{\partial M}$ representing $\sigma$
and the map $\beta$ defined in \S \ref{subsubsec:pont} from
$\Spin(M,s)$ to $\Spin^c\left(M,\beta s\right)$.
\end{remark}

\subsubsection{\spinc--structures as vector fields: the closed case.}

\label{subsubsec:spinc_vs_Euler}
Let $M$ be a closed $3$--manifold. We recall Turaev's definition \cite{TReid}
of an Euler structure on $M$, and how this corresponds to a \spinc--structure on $M$. 

\begin{lemma}
\label{lem:commutes}
The group ${\rm{Spin}^c}(3)$ can be identified with ${\rm{U}}(2)$ in such a way that the diagram
$$
\xymatrix{
*++ {{\rm{SO}}(2)}  \ar@{>->}[d] \ar[r]^{\simeq} & {{\rm{U}}(1)\ } \ar@{>->}[r]
& {{\rm{U}}(2)}\ar[d]^-\simeq \\
{{\rm{SO}}(3)}  & & {{\rm{Spin}^c}(3)} \ar@{->>}[ll]^-{\pi}
}
$$
commutes. Here, $\pi$ is the canonical projection, 
${{\rm{SO}}(2)}$ is identified with ${{\rm{U}}(1)}$ in the usual way and 
is embedded into ${{\rm{SO}}(3)}$ by $A\mapsto (1) \oplus A$,
whereas ${{\rm{U}}(1)}$ is embedded into ${\rm{U}}(2)$ by $A \mapsto A \oplus (1)$.
\end{lemma}
\begin{proof}
There is a well--known way to construct a $2$--fold covering from $\SU(2)$ onto $\SO(3)$,
which consists in identifying $\SU(2)$ with the group of unit quaternions, $\R^3$
with the space of pure quaternions and making the former act on the latter by conjugation. 
Thus, we get a unique group isomorphism $\xymatrix{\SU(2) \ar[r]^-\simeq & \Spin(3)}$ which
makes the two projections onto $\SO(3)$ commute. Then, the isomorphism
$$
\xymatrix{
{\frac{\hbox{\scriptsize SU}(2)\times \hbox{\scriptsize U}(1)}{\Zi_{2}}} \ar[r]^-{\simeq}& {\U(2)}
}
$$
sending $\left[\left(A,z\right)\right]$ to $zA$ induces a group isomorphism 
$\xymatrix{\U(2) \ar[r]^-\simeq & \Spinc(3)}$.
The reader may easily verify the commutativity of  the above diagram.
\end{proof}

\begin{definition}
An \emph{Euler structure} on $M$ is a punctured homotopy class of nonsingular vector fields on $M$. 
Precisely, two nonsingular vector fields $v$ and $v'$ on $M$ are considered as equivalent, 
when there exists a point $x\in M$ such that the restrictions of $v$ and $v'$ to $M\setminus x$ are homotopic
among nonsingular vector fields on $M\setminus x$. 
The set of Euler structures on $M$ is denoted by Eul$(M)$.
\end{definition}
\noindent
If a cellular decomposition of $M$ is given, punctured homotopy
coincides with homotopy on the $2$--skeleton of $M$.
Then, obstruction theory applied to the bundle
of non-zero vectors tangent to $M$ says
that Euler structures do exist (Poincar\'e--Hopf theorem: $\chi(M)=0$) 
and that they form a $H^2\left(M;\pi_2\left(\hbox{T}_yM \setminus 0\right)\right)$--affine space 
(where $y\in M$). Since $M$ has come with an orientation, $\hbox{Eul}(M)$ is naturally a $H^2(M)$--affine space.
\begin{lemma}[Turaev, \cite{TSpinc}]
\label{lem:Eul_spinc}
There exists a canonical $H^2(M)$--equivariant bijection
$$
\xymatrix{
{{\rm{Eul}}(M)} \ar[r]^-{\mu} & {{\rm{Spin}}^c(M)}.
}
$$
\end{lemma}
\begin{proof} 
Let $v$ be a nonsingular vector field on $M$.
We are going to associate to $v$ a \spinc--structure in the usual sense 
(see Lemma \ref{lem:usual_definition}) and, for this, we need to endow $M$ with a metric. 
Orient $\langle v \rangle^{\bot}$, the orthogonal complement of $\langle v \rangle$ in T$M$, 
with the ``right hand'' rule ($v$ being taken as right thumb).  
Then, SO$\left(\langle v \rangle^{\bot}\right)$ is a reduction
of SO$\left(\hbox{T}M\right)$ with respect to the inclusion of SO$(2)$ into SO$(3)$
defined by $A\mapsto (1) \oplus A$. 
The bundle SO$\left(\langle v \rangle^{\bot}\right)$,
together with the homomorphism $\SO(2)\simeq \U(1) \to \U(2)$ 
defined by $A\mapsto A \oplus (1)$,  induces a principal $\U(2)$--bundle $\eta$.
According to Lemma \ref{lem:commutes}, 
$\eta$ can be declared to be a principal $\Spinc(3)$--bundle and can be accompanied
with an isomorphism $H:\eta/\U(1)\to \SO\left(\hbox{T}M\right)$.
The \spinc--structure $\left[(\eta, H)\right]$ on $M$ only depends on the punctured homotopy class of $v$,
and we set $\mu([v])= \left[(\eta, H)\right]$. The map $\mu$ can be verified to be $H^2(M)$--equivariant. 
\end{proof}
\begin{remark}
\label{rem:Chern_Euler}
Let $[v]$ be an Euler structure on $M$.
The isomorphism class of principal $\U(1)$--bundles induced by the \spinc--structure $\mu([v])$ 
in \S \ref{subsubsec:Chern} is represented by $\SO\left(\langle v^{\bot}\rangle\right)$,
since the homomorphism $\Spinc(3) \to \U(1)$ used there corresponds to the determinant map
through the isomorphism $\Spinc(3)\simeq \U(2)$ of Lemma \ref{lem:commutes}.
Consequently, the Chern class of $\mu([v])$ is the Euler class $e(\hbox{T}M/ \langle v \rangle)$,
i.e. the obstruction to find a nonsingular vector field on $M$ transverse to $v$.
\end{remark}

According to the previous remark, \spinc--structures arising from \spin--structures correspond
to nonsingular vector fields on $M$ which can be completed. 

More precisely, let a \emph{parallelization} of $M$  be a punctured homotopy class
of trivializations $t=(t_1,t_2,t_3)$ of the oriented vector bundle T$M$, and denote the 
set of such structures by Parall$(M)$. Obstruction theory applied to the
bundle of oriented frames of $M$ says that parallelizations do exist (Stiefel theorem: $w_2(M)=0$)
and that they form a $H^1(M;\Z_2)$--affine space. (In the case of trivializations of T$M$, 
homotopy on the $2$--skeleton coincides 
with homotopy on the $1$--skeleton since $\pi_2\left(\hbox{GL}_+(3)\right)=0$.) 
Thus, one obtains the following well--known fact \cite{Mi,KiBook}.
\begin{lemma}
\label{lem:parall_spin}
There exists a canonical $H^1(M;\Z_2)$--equivariant bijection
$$
\xymatrix{
{{\rm{Parall}}(M)}\ar[r]^{\mu} & {{\rm{Spin}}(M)}.
}
$$
\end{lemma}
Define a map $\beta: {\hbox{Parall}(M)} \to {\hbox{Eul}(M)}$  by
$\beta([t])=[t_1]$ for any trivialization $t=(t_1,t_2,t_3)$ of T$M$. 
The next lemma follows from the definitions.
\begin{lemma}
\label{lem:vector_beta}
The following diagram is commutative:
$$
\xymatrix{
{{\rm{Parall}}(M)} \ar[r]^{\mu}_\simeq \ar[d]_\beta & {{\rm{Spin}}(M)} \ar[d]^\beta\\
{{\rm{Eul}}(M)}\ar[r]^{\mu}_\simeq & {{\rm{Spin}}^c(M)}.
}
$$
\end{lemma}

\subsubsection{\spinc--structures as vector fields: the boundary case}

\label{subsubsec:Euler_relative}
Let $M$ be a $3$--manifold with boundary.
We define Euler structures on $M$ which are relative to a homotopy class
of trivializations of $\R \oplus \hbox{T}\partial M$. We start with a preliminary observation.

What has been done in \S \ref{subsubsec:spinc_vs_Euler}
for the oriented tangent bundle of a closed $3$--manifold 
works for any $3$--dimensional oriented vector bundle.
In particular, if $S$ is a closed surface, 
\S \ref{subsubsec:spinc_vs_Euler} can be repeated for $\R \oplus \hbox{T}{S}$.
This repetition ends with the following commutative diagram:
$$
\xymatrix{
 {\Parall\left(\R \oplus \hbox{T}S\right)} \ar[r]^{\mu}_\simeq \ar[d]_\beta &
{{\rm{Spin}}(\R \oplus \hbox{T} S)} \ar[d]^\beta
&  \ar[l]_-\simeq {{\rm{Spin}}(S)} \\
{\hbox{Eul}\left( \R \oplus \hbox{T}S\right)} \ar[r]^{\mu}_\simeq &
{{\rm{Spin}}^c(\R \oplus \hbox{T} S)} & \ar[l]_-\simeq  {\Spin^c(S)}.
}
$$
The only change is that, because the base space $S$ is now $2$--dimensional,
homotopies are not punctured anymore.
An \emph{Euler structure} on $\R \oplus \hbox{T}{S}$ is defined to be a homotopy class
of nonsingular sections of this vector bundle and, similarly, a \emph{parallelization}  
on $\R \oplus \hbox{T}{S}$ is a homotopy class of trivializations of this oriented vector bundle. 
\begin{example}
Thus, the section $v=(1,0)$ of $\R \oplus \hbox{T}S$
determines a \spinc--structure $\mu([v])$ on the surface $S$. 
By Remark \ref{rem:Chern_Euler}, the Chern class of $\mu([v])$
coincides with the Euler class $e(\hbox{T}S)$ of the surface $S$.
\end{example}

In the sequel, we fix a parallelization $\tau$ on $\R \oplus \hbox{T}\partial M$.
The observation at the beginning of \S \ref{subsubsec:boundary}
allows us to identify $\tau$ with a homotopy class of trivializations of $\hbox{T}M|_{\partial M}$.

Fix, in this paragraph, a nonsingular section $s$ of T$M|_{\partial M}$.
An \emph{Euler structure on $M$ relative to $s$} is
a punctured homotopy class rel $\partial M$ of nonsingular vector fields on $M$ that extend $s$.
We denote by Eul$(M,s)$ the set of such structures. 
Obstruction theory says that there is an obstruction 
$w(M,s) \in H^3(M,\partial M)$ to the existence of such structures
and, when the latter happens to vanish, that the set Eul$(M,s)$ is naturally
a $H^2(M,\partial M)$--affine space. (Here, again, we use the given orientation of $M$
to make $\Z$ the coefficients group.) As an application of the Poincar\'e--Hopf theorem
and obstruction calculi on the double $M\cup_{\rm{\scriptsize Id}} (-M)$,
one obtains that 
\begin{equation}
\label{eq:obstruction}
2\cdot \langle w(M,s), [M,\partial M]\rangle
= \left\langle e\left(\hbox{T}M|_{\partial M}/ \langle s \rangle\right), [\partial M]\right\rangle \in \Z.
\end{equation}

The following lemma can be proved formally the same way as Lemma \ref{lem:relative_spinc}. 
The first statement is also a direct consequence of (\ref{eq:obstruction}).
\begin{lemma}
\label{lem:relative_Euler}
For any trivialization $t=(t_1,t_2,t_3)$ of ${\rm{T}}M|_{\partial M}$ representing $\tau$
(which we denote by $t \in \tau$), the nonsingular vector field
$t_1$ can be extended to $M$.
Moreover, for any $t,t' \in \tau$, there exists a canonical $H^2(M,\partial M)$--equivariant bijection
$$
\xymatrix{
{{\rm{Eul}}\left(M,t_1\right)} \ar[r]^-{\rho_{t,t'}} & {{\rm{Eul}}(M,t_1')}.
}
$$
Lastly, for any $t,t',t'' \in \tau$, we have that $\rho_{t',t''} \circ \rho_{t,t'} = \rho_{t,t''}$.
\end{lemma}
\begin{definition}
\label{def:relative_Euler}
An \emph{Euler structure on $M$ relative to $\tau$} is a pair $(v,t)$ where $t\in \tau$ 
and $v\in {{\rm{Eul}}\left(M,t_1\right)}$,
two such pairs $(v,t)$ and $(v',t')$ being considered as equivalent when $v'=\rho_{t,t'}(v)$.
The set of such structures is denoted by ${\rm{Eul}}(M,\tau)$ and can naturally be
given the structure of a $H^2(M,\partial M)$--affine space. 
\end{definition}

\begin{remark}
\label{rem:Reeb}
Following Turaev, one can describe concretely
how a  $x\in H^2(M,\partial M)$ acts on a 
$[(v,t)] \in \hbox{Eul}(M,\tau)$. 
Let $P^{-1}x \in H_1(M)$ be represented by a smooth
oriented knot $K \subset \hbox{int}(M)$, and let $v'$ be the vector field
obtained from $v$ by ``Reeb turbulentization'' along $K$
(see \cite[\S 5.2]{TReid}). Then, $(v',t)$ represents $[(v,t)]+x$.
\end{remark}

The following relative version of Lemma \ref{lem:Eul_spinc} can be proved similarly.
\begin{lemma}
\label{lem:relative_Eul_spinc}
There exists a canonical $H^2(M,\partial M)$--equivariant bijection
$$
\xymatrix{
{{\rm{Eul}}(M,\tau)} \ar[r]^-{\mu} & 
{{\rm{Spin}}^c\left(M,\mu(\tau)\right)}.
}
$$
\end{lemma}

\subsubsection{Relative Chern classes}

\label{subsubsec:relative_Chern}
Let $M$ be a $3$--manifold  with boundary and let $\sigma$
be a \spin--structure on $\partial M$. In the relative case too, 
there is a \emph{Chern class} map
$$
\xymatrix{
{{\rm{Spin}}^c(M,\sigma)} \ar[r]^c & {H^2(M,\partial M)}
}
$$
which is affine over the doubling map. It can be defined directly
(using Remark \ref{rem:relative}), or undirectly 
regarding relative \spinc--structures as classes of vector fields (\S \ref{subsubsec:Euler_relative}).
This is done in the next paragraph.

Let $\tau$ be the parallelization on $\R \oplus \hbox{T} \partial M$ 
corresponding to $\sigma$ by $\mu$.
For any trivialization $t$ of $\hbox{T}M|_{\partial M}$ representing $\tau$
and for any nonsingular vector field $v$ on $M$ extending $t_1$,
we can consider the relative Euler class
$$
e\left(\hbox{T} M/ \langle v \rangle, t_2 \right) \in H^2(M,\partial M),
$$
i.e. the obstruction to extend the nonsingular section $t_2$  
of $\hbox{T} M/ \langle v \rangle$ from $\partial M$ to the whole of $M$.
Clearly, this only depends on the equivalence class 
$[(v,t)]$ of $(v,t)$ in the sense of Definition \ref{def:relative_Euler}.
Thus, we get a canonical map
$$
\xymatrix{
{{\rm{Eul}}(M,\tau)} \ar[r] & {H^2(M,\partial M)}
}
$$
which can be verified to be affine over the doubling map thanks to Remark \ref{rem:Reeb}.
Its composition with $\mu^{-1}$ is defined to be $c$.
(Compare with Remark \ref{rem:Chern_Euler}.)

\begin{remark}
\label{rem:relative_Chern}
For any $\alpha \in \Spin^c(M,\sigma)$, the Chern class $c(\alpha)$ vanishes
if and only if $\alpha$ comes from the set $\Spin(M,\sigma)$
defined in Remark \ref{rem:from_spin_to_spinc_relative}.
\end{remark}

We now compute the modulo $2$ reduction of a relative Chern class.
First, recall that the cobordism group $\Omega_1^{\Spinscript}$
is isomorphic to $\Z_2$ \cite{Mi,KiBook}.
For a closed surface $S$, there is the Atiyah--Johnson correspondence
$$
\xymatrix{
{\Spin(S)} \ar[r]^{q}_{\simeq} & {\Quad(S)}
}
$$
between spin structures on $S$ and quadratic functions with the modulo
$2$ intersection pairing of $S$ as associated bilinear pairing \cite{A,J}.
The quadratic function $q_\sigma:H_1(S;\Z_2)\to \Z_2$
corresponding to $\sigma\in \Spin(S)$ is defined by
\begin{displaymath}
q_\sigma\left([\gamma]\right)=\left[\left(\gamma,\sigma|_\gamma\right)\right]
\in \Omega_1^{\Spinscript}\simeq \Z_2
\end{displaymath}
for any oriented simple closed curve $\gamma$ on $S$.
\begin{lemma} 
\label{lem:Chern_mod2}
The following identity holds for any $\alpha\in {\rm{Spin}}^c(M,\sigma)$:
\begin{displaymath}
\forall y\in H_2(M,\partial M), \quad
\langle c(\alpha),y\rangle \textrm{ mod } 2 = q_\sigma\left(\partial_*(y)\right).
\end{displaymath}
Here, $\partial_*:H_2(M,\partial M) \to H_1(\partial M)$
denotes the connecting homomorphism of the pair $(M,\partial M)$ 
and is followed by the modulo $2$ reduction.
\end{lemma}
\begin{proof}
The modulo $2$ reduction of $c(\alpha)$ is
\begin{displaymath} 
w_2(M,\sigma)\in H^2(M,\partial M;\Z_2), 
\end{displaymath}
i.e. the obstruction to extend $\sigma$ to the whole manifold $M$. 
Let $\Sigma$ be a connected immersed surface in $M$ such that 
$\partial \Sigma$ is $\partial M\cap \Sigma$, $\partial \Sigma$ 
has no singularity and $\Sigma$ represents the modulo $2$ reduction of $y$.
Then, $\langle c(\alpha),y\rangle \hbox{ mod } 2= \langle w_2(M,\sigma),[\Sigma]\rangle$ 
is equal to $\langle w_2\left(\Sigma,\sigma|_{\partial \Sigma}\right),
\left[\Sigma\right]\rangle$ and so is the obstruction to extend 
the \spin--structure $\sigma|_{\partial \Sigma}$ to the whole
surface $\Sigma$. Since $\Sigma$ is connected, this is the class 
of $\left(\partial \Sigma, \sigma|_{\partial \Sigma}\right)$
in $\Omega_1^{\Spinscript}$. Thus, we have that 
$\langle c(\alpha),y\rangle \hbox{ mod } 2=
q_\sigma\left(\left[\partial \Sigma\right]\right)= q_\sigma\left(\partial_*(y)\right)$.
\end{proof}
\begin{example} Suppose that $\partial M$ is a disjoint union of tori.
Let $\tau^0$ be the distinguished parallelization corresponding
to the distinguished \spin--structure $\sigma^0$ on the $2$--torus (see Example \ref{ex:tores_spinc}).
An Euler structure on $M$ relative to the union of  copies of
$\tau^0$ is equivalent to a relative Euler structure in the sense of Turaev
\cite[\S 5.1]{TReid}. Lemma \ref{lem:Chern_mod2}
is a generalization of \cite[Lemma 1.3]{TSurg}.
\end{example}

\subsubsection{\spinc--structures as stable complex structures}

We conclude this subsection devoted to the dimension $3$ 
by recalling that, in this case, a \spinc--structure
is equivalent to a stable complex structure on the oriented tangent bundle.
\begin{lemma}
\label{lem:U_dim3}
If $M$ is a closed  $3$--manifold, then the canonical map
$$
\xymatrix{{\rm{U}}^{s}(M) \ar[r]^-\omega &  {{\rm{Spin}}^c}(M)}
$$
introduced in \S \ref{subsubsec:U} is bijective.
\end{lemma} 
\begin{proof}
Endow $M$ with a Riemannian metric and 
consider a nonsingular vector field $v$ on $M$. Then,
$\R\oplus \hbox{T}M$ splits as $\left(\R\oplus \langle v \rangle\right)\oplus \langle v \rangle^\bot$,
which is the sum of two oriented $2$--dimensional vector bundles.
So, via the inclusion of $\U(1)\times \U(1)$ into $\U(2)$ defined by $(A,B) \mapsto (A) \oplus (B)$,
$v$ defines a complex structure $J_v$ on $\R \oplus \hbox{T}M$. Thus, we get a map
from Eul$(M)$ to the set of stable complex structures on T$M$ up to punctured homotopy. 
By obstruction theory applied to the fibration $\hbox{B} \U \to \hbox{B}\SO$ 
with fiber type $\SO/\U$, the latter set is a $H^2(M)$--affine space
and that map is $H^2(M)$--equivariant. Thus, since $\pi_3(\SO/\U)$ is zero, we get a bijective map
$$
\xymatrix{\hbox{Eul}(M) \ar[r]^-J_-\simeq & \hbox{U}^{s}(M)}.
$$
It can be verified that $\omega \circ J$ is the map $\mu$ from Lemma \ref{lem:Eul_spinc}. 
(This verification amounts to checking that some two group homomorphisms from
$\U(1)$ to $\Spinc(4)$ coincide.)
\end{proof}

\subsection{Gluing of complex spin structures}

\label{subsec:gluing}
In this subsection, we deal with the technical problem of gluing \spinc--structures.
We formulate the gluing in terms of (rigid) \spinc--structures, 
but the reader may easily translate the statement and the proof
in terms of vector fields and Euler structures.
 
Let $M$ be a closed $n$--manifold obtained by
gluing two $n$--manifolds $M_1$ and $M_2$ along their boundaries:
\begin{displaymath}
M=M_1 \ \cup_f \ M_2.
\end{displaymath}
This involves a positive diffeomorphism 
$f:-\partial M_{2} \to \partial M_1$ as well as 
a collar neighborhood of $\partial M_i$ in $M_i$.
The inclusion $M_i\hookrightarrow M$ will be denoted by $j_i$.
\begin{lemma}
\label{lem:general_gluing}
For $i=1,2$, let $s_i$ be a rigid \spinc--structure on ${\rm{T}} M_i |_{\partial M_{i}}$.
Having identified $\R \oplus  {\rm{T}} \partial M_{i}$ with ${\rm{T}} M_i |_{\partial M_{i}}$
thanks to the collar, we assume that $s_1 \circ \left( -{\rm{Id}} \oplus {\rm{T}}f \right) =s_2$. 
If the relative obstructions $w(M_i,s_i)$'s vanish,
then  the absolute obstruction $w(M)$ does too and there is a canonical \emph{gluing} map 
\begin{displaymath}
\xymatrix{
{{\rm{Spin}}^c(M_1,s_1) \times {\rm{Spin}}^c(M_2,s_2)} 
\ar[r]^-{\cup_f} & {\rm{Spin}}^c(M)
}
\end{displaymath}
which is affine over
\begin{displaymath}
\xymatrix{
{H^2(M_1,\partial M_1)\oplus H^2(M_2,\partial M_2)}\ar@{.>}[r] 
\ar[d]_{P^{-1} \times P^{-1}} & {H^2(M)} \\
{H_{n-2}(M_{1})\oplus H_{n-2}(M_{2})}
\ar[r]^-{j_{1,*}\oplus j_{2,*}} & H_{n-2}(M)\ar[u]_P.
}
\end{displaymath}
\end{lemma}
\begin{proof}
For $i=1,2$, let $\alpha_i\in \Spin^c(M_i,s_i)$ be represented
by a rigid structure $a_i$. The structures $a_1$ and $a_2$ can be glued together
by means of T$f$: we obtain a rigid \spinc--structure
on $M$ whose homotopy class does not depend on the choices of
$a_1$ and $a_2$ in $\alpha_1$ and $\alpha_2$ respectively. 
We denote it by $\alpha_1\cup_f \alpha_2 \in \Spin^c(M)$.

Let us prove that this map $\cup_f$ is affine. For $i=1,2$,
let $\mathcal{C}_i$ be a smooth triangulation of $M_i$ such that
$\mathcal{C}_1|_{\partial M_1}$ corresponds to 
$\mathcal{C}_2|_{\partial M_2}$ by $f$.
We denote by $\mathcal{C}_{i}^*$ the cellular decomposition of $M_i$
dual to the triangulation $\mathcal{C}_{i}$.

On the one hand, we consider the \emph{union} $\mathcal{C}$ of the 
triangulations $\mathcal{C}_1$ and $\mathcal{C}_2$: 
a simplex of $\mathcal{C}$ 
is a simplex of $\mathcal{C}_i$ for $i=1$ or $2$, and 
simplices of $\partial M_1$ are identified with simplices of
$\partial M_2$ by $f$. On the other hand, we consider 
the \emph{gluing} $\mathcal{C}^*$ 
of the cellular decompositions $\mathcal{C}^*_1$ and $\mathcal{C}^*_2$: 
a cell of $\mathcal{C}^*$ either is a cell of $\mathcal{C}^*_i$
which does not intersect $\partial M_i$, either is the gluing by $f$ 
of a cell belonging to $\mathcal{C}^*_1$ with a cell of $\mathcal{C}^*_2$ 
along a face lying in $\partial M_1\cong-\partial M_2$.
Then, $\mathcal{C}$ is a smooth triangulation of $M$
and $\mathcal{C}^*$ is its dual cellular decomposition.
Cohomology will be calculated with $\mathcal{C}$
while homology will be computed with $\mathcal{C}^*$.

For $i=1,2$, consider some $\alpha_i , \alpha'_i \in \Spin^c(M_i,s_i)$ 
and set $\alpha=\alpha_{1}\cup_f \alpha_2$ and 
$\alpha'=\alpha'_1\cup_f \alpha'_2$. We want to prove the following equality:
\begin{equation}
j_{1,*}P^{-1}\left(\alpha_1-\alpha'_1\right)+
j_{2,*}P^{-1}\left(\alpha_2-\alpha_2'\right)
= P^{-1}\left( \alpha-\alpha'\right)\in H_{n-2}(M).
\end{equation}
For $i=1,2$, let $a_i,a'_i\in \Spin^c_r(M_i)$ represent $\alpha_i$ and $\alpha_i'$ respectively
and coincide on the $1$--skeleton of $\mathcal{C}_i$ (and, of course, on $\partial M_i$). 
Suppose that we have fixed a morphism of oriented vector bundles
$\hbox{T}M_i \to \gamma_{\SOscript(n)}$:
then, the rigid structures $a_i$ and $a'_i$ can be identified 
with lifts $M_i \to \hbox{B}\Spin^c(n)$  by B$\pi$ of the base maps $M_i \to \hbox{B}\SO(n)$. 
The obstruction $\alpha_i-\alpha'_i\in H^2(M_i,\partial M_i)$ 
is the class of the $2$--cocycle which assigns
to each $2$--simplex $e^{i}_k$ of $\mathcal{C}_i$ 
outside $\partial M_i$,
this element $z^i_k$ of $\pi_2(\hbox{B}\U(1))\simeq \pi_2(K(\Z,2))\simeq \Z$ 
obtained by gluing $a_{i}|_{e_k^{i}}$ and $a'_{i}|_{e_k^{i}}$ 
along $\partial e_k^{i}$.
So, we  have that $P^{-1}(\alpha_i - \alpha'_i)$
$=\left[\sum_{k} z_k^i \cdot e_k^{*,i}\right]$
if $e_k^{*,i}$ denotes the $(n-2)$--cell dual to $e_k^{i}$.

Moreover, $a:=a_1\cup_f a_2$ and $a':=a'_1\cup_f a'_2$ represent $\alpha$ and $\alpha'$ respectively.
Using these rigid structures, we can describe explicitely a
$2$--cocycle representing $\alpha-\alpha'$ as well. This $2$--cocycle
sends any $2$--simplex of $\mathcal{C}_1\cup_f\mathcal{C}_2$ contained
in $\partial M_1\cong-\partial M_2$ to $0\in\Z$ so that 
$P^{-1}(\alpha-\alpha')$ is represented by 
$\sum_k z_k^1 \cdot e_k^{*,1} + \sum_k z_k^2 \cdot e_k^{*,2}$.
\end{proof}

Suppose now that the manifolds have dimension $n=3$.
This is the gluing lemma that we will use in the next sections.
\begin{lemma}
\label{lem:gluing}
Let $\sigma_1 \in {\rm{Spin}}(\partial M_1)$ and 
$\sigma_2 \in {\rm{Spin}}(\partial M_2)$ 
be such that $f^{*}(\sigma_1)=-\sigma_2$. 
Then, there is a canonical \emph{gluing} map
$$
\xymatrix{
{{\rm{Spin}}^c(M_1,\sigma_1) \times {\rm{Spin}}^c(M_2,\sigma_2)}
\ar[r]^-{\cup_f} & {\rm{Spin}}^c(M)
}
$$
which is affine over
$$
\xymatrix{
{H^2(M_1,\partial M_1)\oplus H^2(M_2,\partial M_2)}\ar@{.>}[r]
\ar[d]_{P^{-1} \times P^{-1}} & {H^2(M)} \\
{H_{1}(M_{1})\oplus H_{1}(M_{2})}
\ar[r]^-{j_{1,*}\oplus j_{2,*}} & H_{1}(M). \ar[u]_P
}
$$
Moreover, for any $\alpha_1 \in {{\rm{Spin}}}^c(M_1,\sigma_1)$ and 
$\alpha_2 \in {{\rm{Spin}}}^c(M_2,\sigma_2)$, the following identity between Chern classes holds:
\begin{displaymath}
P^{-1}c(\alpha_1\cup_f \alpha_2)= 
j_{1,*}P^{-1}c(\alpha_1)+j_{2,*}P^{-1}c(\alpha_2) 
\in H_1(M).
\end{displaymath}
\end{lemma}
\begin{proof}
Choose a rigid \spin--structure $s_1$ on $\hbox{T}M_1|_{\partial M_1}$ representing
$\sigma_1$, which we denoted by $s_1\in \sigma_1$. This induces a $s_2\in \sigma_2$
by setting $s_2=s_1 \circ \left( -{\rm{Id}} \oplus {\rm{T}}f \right)$.
By Lemma \ref{lem:relative_spinc}, the obstructions $w(M_i,\beta s_i)$'s vanish and so, 
by Lemma \ref{lem:general_gluing}, there is a gluing map 
with domain ${{\rm{Spin}}^c(M_1,\beta s_1) \times {\rm{Spin}}^c(M_2,\beta s_2)}$.

Another choice $s'_1\in \sigma_1$ would induce another $s'_2\in \sigma_2$ 
and would lead to another gluing map this time
with domain ${{\rm{Spin}}^c(M_1,\beta s_1') \times {\rm{Spin}}^c(M_2,\beta s_2')}$.
Nevertheless, using the ``double collar'' of $\partial M_1 \cong -\partial M_2$ in $M$,
one easily sees that the identifications $\rho_{s_1,s'_1}$ and $\rho_{s_2,s'_2}$ 
from Lemma \ref{lem:relative_spinc} make those two gluing maps agree. 

The first assertion of the lemma then follows. The second one is proved 
with arguments similar to those used in the proof of Lemma \ref{lem:general_gluing}
(gluing of obstructions in compact oriented manifolds using Poincar\'e duality).
\end{proof}
\begin{remark}
\label{rem:boundary}
If $M$ is obtained by gluing $M_1$ and $M_2$ along only
part of their boundaries (so that $\partial M\neq \varnothing$),
Lemma \ref{lem:gluing} can easily be generalized
to produce \spinc--structures on $M$ relative to a fixed \spin--structure on its boundary.
\end{remark}

\section{Linking quadratic function of a three--manifold with complex spin structure} 

\label{sec:quad_spinc}
In this section, we define the quadratic function $\phi_{M,\sigma}$
associated to a closed oriented $3$--manifold $M$ equipped with a \spinc--structure $\sigma$.
We present its elementary properties and connect it to previously known constructions.

\subsection{Quadratic functions on torsion Abelian groups}

\label{subsec:discriminant}
We fix some notations. If $A$ and $B$ are Abelian groups and 
if $b:A\times A \to B$ is a symmetric bilinear pairing, 
we denote by $\widehat{b} : A \to {\textrm{Hom}}(A,B)$ the adjoint map.
The pairing $b$ is said to be \emph{nondegenerate} (respectively \emph{nonsingular})
if $\widehat{b}$ is injective (respectively bijective).
We denote by $A^*$ the group $\hbox{Hom}(A,\Z)$ when $A$ is free, 
the group $\hbox{Hom}(A,\Q)$ when $A$ is a $\Q$--vector space and
the group $\hbox{Hom}(A,\Q/\Z)$ when $A$ is torsion.
Lastly, application of the functor $-\otimes \Q$ is indicated by a subscript $\Qi$.

\subsubsection{Basic notions about quadratic functions}

Let  $G$ be a torsion Abelian group. 

A map $q:G\to \Q/\Z$ is said to be a \emph{quadratic function} on $G$ if 
$$
b_q(x,y)=q(x+y)-q(x)-q(y)
$$ 
defines a (symmetric) bilinear pairing $b_q:G\times G\to \Q/\Z$. 
The quadratic function $q$ is said to be \emph{nondegenerate} 
if $b_q$ is nondegenerate, and \emph{homogeneous} if $q(-x)=q(x)$ for any $x\in G$.
Apart from the bilinear pairing $b_q$, one can associate to $q$ its \emph{radical}
$$\hbox{Ker}(q)=\hbox{Ker}\ \widehat{b_q}\subset G,$$ 
its \emph{homogeneity defect}
$$d_q: G\to \Q/\Z, \ x\mapsto q(x)-q(-x)$$ 
and, in case when $G$ happens to be finite, its \emph{Gauss sum} 
$$\gamma(q)=\sum_{x\in G} \exp\left({2i\pi q(x)}\right) \in \C.$$

Given a symmetric bilinear pairing $b:G\times G\to \Q/\Z$, we say
that $q:G\to \Q/\Z$ is a quadratic function \emph{over $b$} if $b_q=b$. 
The group $G^*$ acts freely and transitively on $\Quad(b)$, the set of quadratic functions 
over $b$, just as maps $G\to \Q/\Z$ add up. So, $\Quad(b)$ is a $G^*$--affine space.

There is a procedure to produce quadratic functions on torsion Abelian groups,
known as the ``discriminant'' construction. 

\subsubsection{The discriminant construction}

In the litterature, the discriminant construction is usually restricted to nondegenerate 
bilinear lattices and produces quadratic functions on finite Abelian groups.
The general case has been considered in \cite{DM1}, to which we refer for details and proofs. 
Here, we briefly review the construction.\\

A \emph{lattice} $H$ is a free finitely generated Abelian group.
A \emph{bilinear lattice} $(H,f)$ is a symmetric bilinear pairing
$f:H\times H\to \Z$ on a lattice $H$. Let also 
$$H^\sharp=\{x\in H_{\Qi} : f_\mathbf{\Qi}(x,H)\subset \Z\}$$ 
be the dual lattice. A \emph{Wu class} for $(H,f)$ is an element $w\in H$ such that
$$\forall x\in H, \ f(x,x)-f(w,x) \in 2\Z.$$
A \emph{characteristic form} for $(H,f)$ is an element $c\in H^*=\Hom(H,\Z)$ satisfying
$$\forall x\in H, \ f(x,x)-c(x) \in 2\Z.$$
The sets of characteristic forms and Wu classes for $(H,f)$ are denoted by
$\hbox{Char}(f)$ and $\hbox{Wu}(f)$ respectively. Those sets are not empty and are related
by the map $w\mapsto \widehat{f}(w),\hbox{Wu}(f) \to \hbox{Char}(f)$.\\

Let $(H,f)$ be a bilinear lattice. Consider the torsion Abelian group $G_f=H^\sharp/H$ and the map 
$$L_f: G_f\times G_f \to \Q/\Z, \ \left([x],[y]\right) \mapsto f_{\Qi}(x,y) \textrm{ mod } 1.$$
The pairing $L_f$ is symmetric and bilinear, with radical
$\hbox{Ker}\ \widehat{L_f} \simeq \left( \hbox{Ker} \widehat{f}\right) \otimes \Q/\Z$.

Observe that the adjoint map $\widehat{f_{\Qi}}:H_{\Qi} \to H_{\Qi}^*$ 
restricted to $H^{\sharp}$ induces an epimorphism $G_f \to \hbox{Tors}\ {\hbox{Coker}}\ \widehat{f}$.
Hence the short exact sequence
\begin{equation}
0 \to  \hbox{Ker}\ \widehat{L_f} \to  
G_f \to \hbox{Tors}\ {\hbox{Coker}}\ \widehat{f} \to 0, 
\label{eq:shortexact}
\end{equation}
which can be verified to split (non-canonically). Therefore, $G_f$ is the direct sum
of a finite Abelian group with as many copies of $\Q/\Z$ 
as the rank of $\hbox{Ker}\ \widehat{f}$.
It follows also from (\ref{eq:shortexact}) that the pairing $L_f$ factors 
to a nondegenerate symmetric bilinear pairing
$$
\xymatrix{{\hbox{Tors}}\ {\hbox{Coker}}\ \widehat{f} 
\times {\hbox{Tors}}\ {\hbox{Coker}}\ \widehat{f} \ar[r]^-{\lambda_f} &\Q/\Z.}
$$

The bilinear map  $H^{*} \times H^{\sharp} \to \Q$ defined 
by $(\alpha, x) \mapsto \alpha_{\Qi}(x)$ induces a bilinear pairing
\begin{equation}
\label{eq:evaluation}
\xymatrix{\hbox{Coker}\ \widehat{f} \times G_{f} \ar[rr]^-{\langle -, -\rangle} && \Q/\Z}
\end{equation} 
which is left nondegenerate and right nonsingular. 
It is left nonsingular if and only if $f$ is nondegenerate.\\

Let now $(H,f,c)$ be a bilinear lattice equipped with a 
characteristic form $c \in H^{*}$. One can associate to this triple
a quadratic function over $L_f$, namely
$$
\phi_{f,c}: G_f \to \Q/\Z, \
[x] \mapsto \frac{1}{2}(f_{\Qi}(x,x)-c_{\Qi}(x))\ \hbox{mod}\ 1.
$$
\begin{definition}
The assignation $(H,f,c)\mapsto (G_f,\phi_{f,c})$ is called 
the \emph{discriminant} construction.
\end{definition}
Let us make a few observations about this construction. First,
note that $\phi_{f,c}$ depends on $c$ only mod $2\widehat{f}(H)$. 
Second, the Abelian group  $H^*/\widehat{f}(H) = 
\hbox{Coker}\ \widehat{f}$ acts freely and transitively on 
Char$(f)/2\widehat{f}(H)$ by setting
$$
\forall [\alpha] \in \hbox{Coker}\ \widehat{f},\ \forall [c] \in \hbox{Char}(f)/2\widehat{f}(H), \
[c] + [\alpha] = [c + 2\alpha] \in \hbox{Char}(f)/2\widehat{f}(H).
$$
Third, since $\hbox{Ker}\ \widehat{L_f}$ is canonically isomorphic to 
$\left(\hbox{Ker}\ \widehat{f}\right)\otimes \Q/\Z$, 
any form Ker  $\widehat{f} \to \Z$ induces a homomorphism 
$\hbox{Ker}\ \widehat{L_f} \to \Q/\Z$. Thus, we get a homomorphism
$j_f:\left(\hbox{Ker}\ \widehat{f}\right)^{*} \to \left(\hbox{Ker}\ \widehat{L_f}\right)^{*}$.
\begin{theorem} \cite{DM1} 
\label{th:algebraic-embedding}
The assignation $c \mapsto \phi_{f,c}$ induces an embedding 
$$
\xymatrix{ {\rm{Char}}(f)/2\widehat{f}(H)\ \ar@{^{(}->}[r]^-{\phi_f} &{\rm{Quad}}(L_{f})}
$$
which is affine over the opposite of the left adjoint of the pairing (\ref{eq:evaluation}).
Moreover, a function $q \in {\rm{Quad}}(L_{f})$ belongs to ${\rm{Im}}\ \phi_{f}$ if and only if 
$q|_{{\rm{\scriptsize Ker}}\ \widehat{L_{f}}}$ belongs to ${\rm{Im}}\ j_{f}$.
\end{theorem}
\begin{remark}
The map $\phi_f$ is bijective if and only if $f$ is nondegenerate.
\end{remark}

We now use the algebraic notions above
as combinatorial descriptions of topological notions.

\subsection{Combinatorial descriptions associated to a surgery presentation} 

\label{subsec:combinatorial}
In this subsection, we fix an ordered oriented framed $n$--component link $L$ in $\mathbf{S}^{3}$. 

We call $V_{L}$ the $3$--manifold obtained from $\mathbf{S}^{3}$
by surgery along $L$ and we denote by $W_{L}$ the 
\emph{trace} of the surgery:
$$
V_L=\partial W_L \quad \hbox{with} \quad
W_L=\mathbf{D}^4\cup \bigcup_{i=1}^n
\left(\mathbf{D}^2\times \mathbf{D}^2\right)_i
$$
where the $2$--handle $\left(\mathbf{D}^2\times \mathbf{D}^2\right)_i$
is attached by embedding $-\left(\mathbf{S}^1\times 
\mathbf{D}^2\right)_i$ into $\mathbf{S}^3=\partial \mathbf{D}^4$
in accordance with the specified framing and orientation of $L_i$.

\begin{figure}[h]
\begin{center}
\includegraphics[width=8cm,height=6cm]{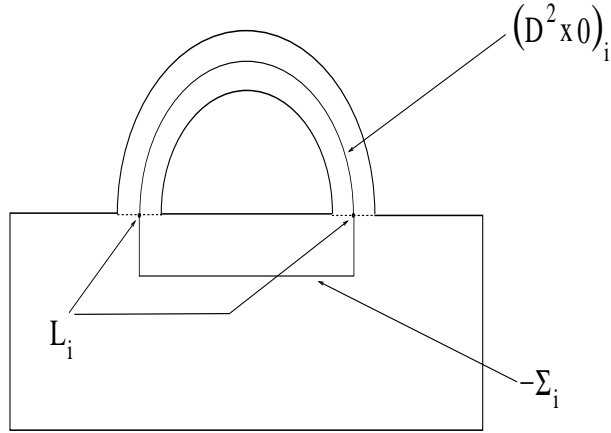}
\caption{The preferred basis of $H_2\left(W_L\right)$.}
\label{fig:handle}
\end{center}
\end{figure}

The group $H_{2}(W_{L})$ is free Abelian of rank $n$, and is given 
the \emph{preferred} basis $([S_{1}],\dots,[S_{n}])$ defined as follows. 
The closed surface $S_i$ is taken to be 
$\left(\mathbf{D}^2\times 0\right)_i
\cup \left(-\Sigma_i\right)$, where $\Sigma_i$ is a Seifert surface
for $L_i$ in $\mathbf{S}^3$ which has been pushed off into
the interior of $\mathbf{D}^4$ as shown in Figure \ref{fig:handle}.
The group $H^2(W_L)$ is identified with
${\Hom}(H_2(W_L),\Z)$ by Kronecker evaluation, and is given the dual basis. In the sequel, we simplify the
notations by setting $H=H_2(W_L)$ (so that $H^2(W_L)$ is identified with $H^*$) 
and by denoting by $f: H\times H \to\Z$ the intersection pairing of $W_L$.
The matrix of $f$ relatively to the preferred basis of $H$ is the linking matrix
$$
B_{L}=\left(b_{ij}\right)_{i,j=1}^n
$$ 
of $L$. Since $(H,f)$ is a bilinear lattice, the constructions of \S \ref{subsec:discriminant} apply.

\subsubsection{Combinatorial description of \spin--structures}

We recall a combinatorial description of $\hbox{Spin}(V_L)$ due to Blanchet \cite{Blanchet}.
Define the set
\begin{displaymath}
\mathcal{S}_L=
\left\{[r]=\left([r_i]\right)_{i=1}^n\in (\Z_2)^n \ : \
\forall i=1,\dots,n,\ \sum_{j=1}^nb_{ij}r_j\equiv b_{ii} 
\textrm{ mod } 2\right\}.
\end{displaymath}
The elements of $\mathcal{S}_L$ are called 
\emph{characteristic solutions} of $B_L$. 
\begin{lemma}
\label{lem:combinatorial_spin}
There are canonical bijections
$$
\xymatrix{
{{\rm{Spin}}(V_L)} \ar[r]_-{\simeq} & 
{{\rm{Wu}}(f)/2H} \ar[r]_-{\simeq} &{\mathcal{S}_L}.
}
$$
\end{lemma}
\noindent
Thus, $\mathcal{S}_L$ shall be refered to as the
\emph{combinatorial description of} ${\rm{Spin}}(V_L)$.
A refined Kirby's theorem dealing with surgery presentations of closed
$3$--dimensional \spin--manifolds can be derived 
from this lemma \cite[Theorem (I.1)]{Blanchet}.
\begin{proof}[Proof of Lemma \ref{lem:combinatorial_spin}]
The preferred basis of $H$ induces an isomorphism
$H/2H\simeq (\Z_2)^n$: the bijection between 
$\hbox{Wu}(f)/2H$ and $\mathcal{S}_L$ is obtained this way.
We now describe a bijection between $\Spin(V_L)$ and $\hbox{Wu}(f)/2H$.
Let $\sigma$ be a \spin--structure on $V_L$.
The obstruction $w_2(W_L,\sigma)$ to extend $\sigma$ to $W_L$
belongs to the group $H^2(W_L,V_L;\Z_2)\simeq H_2(W_L;\Z_2)\simeq H/2H$.
Since $w_2(W_L,\sigma)$ is sent to $w_2(W_L)$ by the restriction map 
$H^2(W_L,V_L;\Z_2)\rightarrow H^2(W_L;\Z_2)$, a representative for
$w_2(W_L,\sigma)$ in $H$ has to be a Wu class for $f$. 
\end{proof}

\subsubsection{Combinatorial description of \spinc--structures}

\label{subsubsec:combinatorial_spinc}
Define the set
\begin{displaymath}
\mathcal{V}_L=
\frac{\left\{ s=(s_i)_{i=1}^n\in\Z^n \ : \ 
\forall i=1,\dots,n,\ s_i \equiv b_{ii}\textrm{ mod }2\right\}}
{2\cdot {\rm{Im}}\ B_L},
\end{displaymath} 
the elements of which are called \emph{Chern vectors} of $B_L$. According to the following lemma,
this set shall be referred to as the \emph{combinatorial description of} Spin$^c(V_L)$.
\begin{lemma}
\label{lem:combinatorial_spinc}
There are canonical bijections
$$
\xymatrix{{{\rm{Spin}}^c(V_L)} \ar[r]_-{\simeq} & 
{{\rm{Char}}(f)/2 \widehat{f}(H)} 
\ar[r]_-{\simeq} &{\mathcal{V}_L}.}
$$
\end{lemma}
\begin{proof}
The preferred basis of $H$ defines an isomorphism $H^*\simeq \Z^n$,
which induces a bijection between $\hbox{Char}(f)/2\widehat{f}(H)$ 
and $\mathcal{V}_L$.
The restriction map ${\Spin}^c(W_L)$
$\rightarrow \Spin^c(V_L)$ is affine
over the map $H^2(W_L)\rightarrow H^2(V_L)$ induced by inclusion.
By exactness of the pair $(W_L,V_L)$, the latter 
is surjective and its kernel coincides with the image
of $\widehat{f}: H\rightarrow H^*$ (by Poincar\'e duality). 
Moreover, since $H^2(W_L)$ is free Abelian, 
a \spinc--structure  on $W_L$ is determined by its Chern class in $H^2(W_L)\simeq H^*$.
Such a class has to be a characteristic form for $f$ 
since its modulo $2$ reduction coincides
with the second Stiefel--Whitney class $w_{2}(W_{L})\in
H^2(W_L;\Z_2)\simeq \hbox{Hom}(H,\Z_2)$. 
Therefore, there is a bijection between ${\Spin}^c(V_L)$ and 
$\hbox{Char}(f)/2\widehat{f}(H)$ defined by
$\sigma\mapsto \left[c\left(\tilde{\sigma}\right)\right]$ where $\tilde{\sigma}$
is an extension of $\sigma$ to $W_L$. (This extension exists since
$w(W_L,\sigma)$ lives in $H^3(W_L,V_L)=0$, see Proposition \ref{prop:relative_obstruction}.)
\end{proof}

If the Chern vector $[s]$ corresponds to the \spinc--structure $\sigma$,
we say that $(L,[s])$ is a \emph{surgery presentation} of the closed
$3$--dimensional \spinc--manifold  $(V_L,\sigma)$. 
On a diagram, we draw the framed link $L$ using the blackboard framing convention,
indicate its orientation and decorate each of its components $L_i$ with the integer $s_i$.

Next, Kirby's theorem \cite{Ki} can easily be extended to deal with 
surgery presentations of \spinc--manifolds. This Spin$^c$ version of Kirby's calculus
will be used in the next section.

\begin{theorem}
\label{th:Spin^c_Kirby}
Let $L$ and $L'$ be ordered oriented framed links in $\mathbf{S}^{3}$.
Equip them with Chern vectors $[s]$ and $[s']$, which
correspond to \spinc--structures $\sigma$ and $\sigma'$ on $V_L$ and $V_{L'}$ respectively. 
Then, the \spinc--manifolds $(V_L,\sigma)$ and $(V_{L'},\sigma')$ 
are \spinc--diffeomorphic if and only if the pairs $(L,[s])$ and $(L',[s'])$ are,
up to re-ordering and up to isotopy,
related one to the other by a finite sequence of  the moves drawn on Figure \ref{fig:Kirby}.
\begin{figure}[h!]
\begin{center}
\includegraphics[width=10cm,height=8cm]{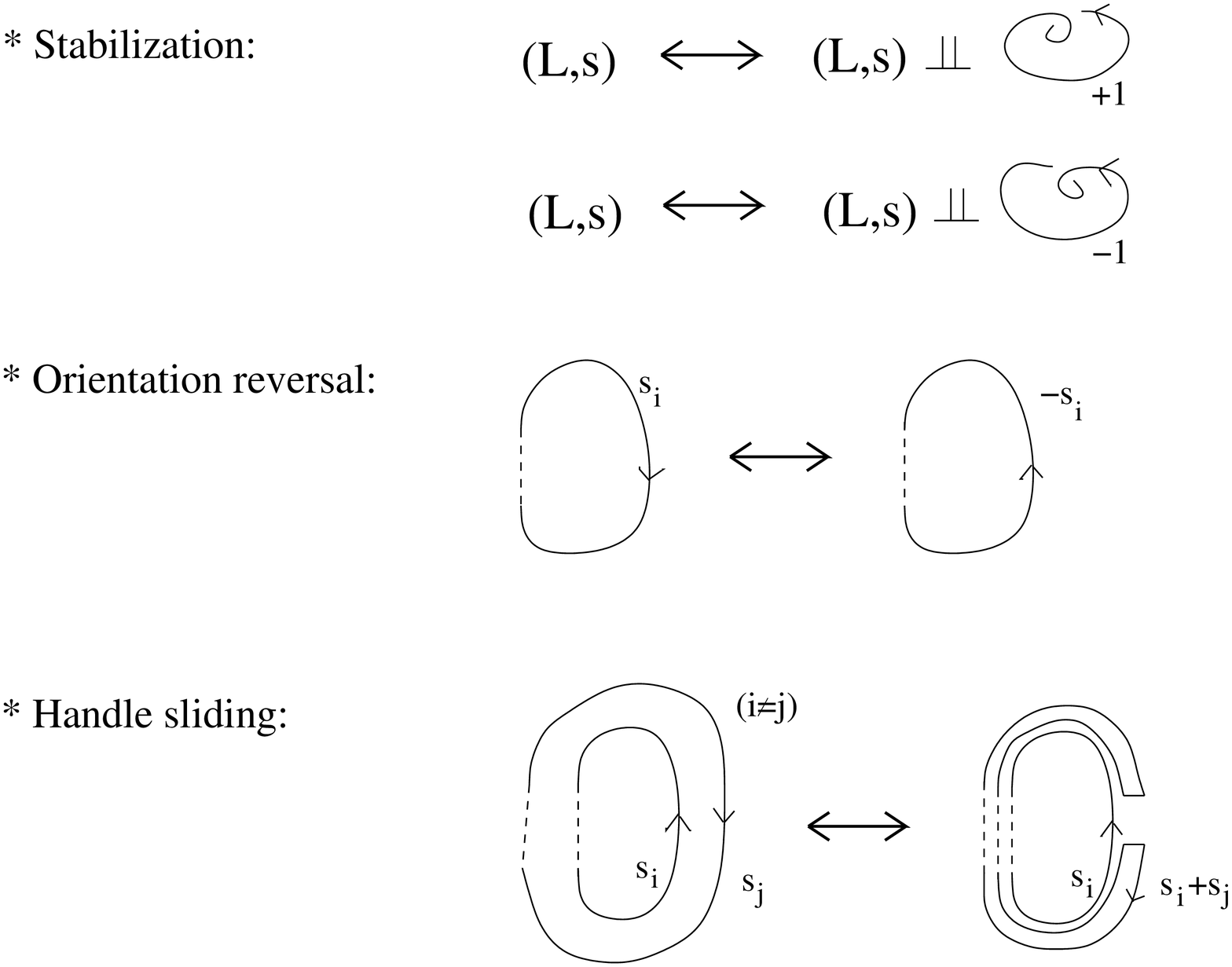}
\caption{$\Spin^c$ Kirby's moves. (Recall that the blackboard framing convention is used, and
that labels refer to Chern vectors.)}
\label{fig:Kirby}
\end{center}
\end{figure}
\end{theorem}

\begin{proof}
This follows from the usual Kirby's theorem. It suffices to show that,
for each Kirby's move $L_1 \to L_2$, the corresponding canonical 
diffeomorphism $V_{L_1}\rightarrow V_{L_2}$ acts 
at the level of \spinc--structures as
combinatorially described on Figure \ref{fig:Kirby}. 
This is a straightforward verification.
\end{proof}
\begin{example}
\label{ex:slam_dunk}
Look at the \emph{slam dunk} move depicted on Figure \ref{fig:slamdunk}.
Here, we are considering the ordered union $L\cup \left(K_1,K_2\right)$
of a $n$--component ordered oriented framed link $L$ 
with an oriented framed knot $K_1$ 
together with its oriented meridian $K_2$. The move is
$$
\left( L\cup \left(K_1,K_2\right), \left[(s_1,\dots,s_n,y,0)\right]\right)
\longleftrightarrow
\left(L,\left[(s_1,\dots,s_n)\right]\right),
$$
where $y$ is the framing number of $K_1$.
It relates two closed \spinc--manifolds which are \spinc--diffeomorphic, as can be shown by re-writing 
the proof of \cite[Lemma 5]{FR} with \spinc{} Kirby's calculi.
\begin{figure}[!h]
\begin{center}
\includegraphics[width=8cm,height=4cm]{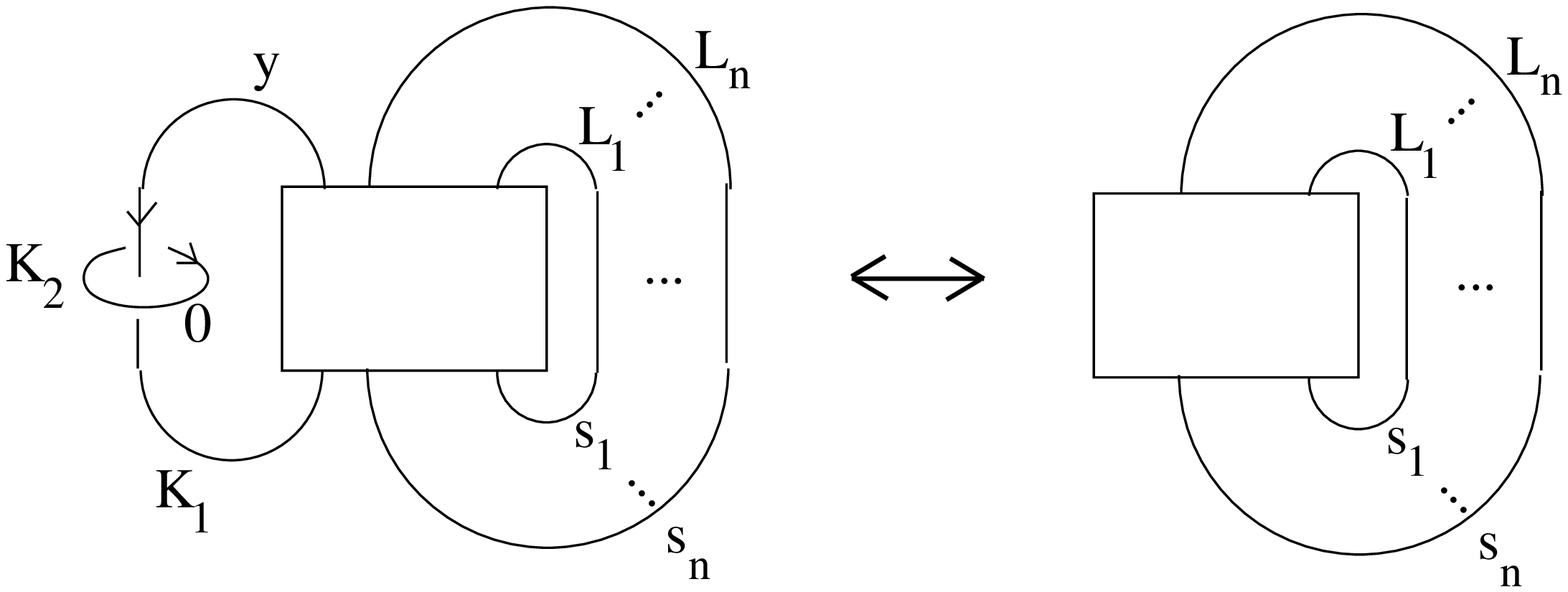}
\caption{\spinc{} slam dunk move.}
\label{fig:slamdunk}
\end{center}
\end{figure}
\end{example}
\begin{remark}
\label{rem:combinatorial_Chern}
There exists a canonical isomorphism $\varrho: \hbox{Coker}\ \widehat{f} \rightarrow H^2(V_L)$,
as defined by the following commutative diagram:
\begin{displaymath}
\xymatrix{
{H^2(W_L,V_L)} \ar[r] & {H^2(W_L)} \ar[r] & H^2(V_L) \ar[r] & 0 \\
{H}\ar[u]_P^\simeq \ar[r]^{\widehat{f}} & H^* \ar[u]^\simeq\ar[r] & 
{\hbox{Coker}\ \widehat{f}}\ar@{.>}[u]^\simeq_\varrho \ar[r]& 0.
}
\end{displaymath}
Then, the affine action of $H^2(V_L)$ on ${\Spin}^c(V_L)$ writes combinatorially:
$$
\forall [x] \in \hbox{Coker}\ \widehat{f}, \
\forall [c] \in \hbox{Char}(f)/2\widehat{f}(H), \
[c]+[x] = [c+2x].
$$
The Chern class map $c:{\Spin}^c(V_L)\rightarrow H^2(V_L)$
is combinatorially described by the map 
$c:\hbox{Char}(f)/2\widehat{f}(H)
\rightarrow\hbox{Coker}\ \widehat{f},\ [c]\mapsto [c]$.
\end{remark}

\subsubsection{From ${\rm{Spin}}$ to ${\rm{Spin}}^c$ in a combinatorial way}

\label{subsubsec:pont_combinatorial}
We now relate the combinatorial description of $\hbox{Spin}(V_L)$
to that of $\hbox{Spin}^c(V_L)$.
\begin{lemma}
The canonical map $\beta:{\rm{Spin}}(V_L)\rightarrow {\rm{Spin}}^c(V_L)$
corresponds to the map 
$\beta:{\rm{Wu}}(f)/2H \rightarrow {\rm{Char}}(f)/2\widehat{f}(H)$ defined by 
$\beta([w])=\left[\widehat{f}(w)\right]$
or, equivalently, to the map $\beta:{\mathcal{S}}_L\rightarrow
{\mathcal{V}_L}$ defined by $\beta\left([r]\right)=
\left[B_L\cdot r\right]$.
\end{lemma}
\begin{proof}
Take $\sigma\in \hbox{Spin}(V_L)$ and let 
$r_\sigma \in H^2(W_L,V_L)\simeq\Z^n$ be an integral representative 
for the obstruction $w_2\left(W_L,\sigma\right) \in H^2(W_L,V_L;\Z_2)
\simeq \left(\Z_2\right)^n$ to extend $\sigma$ to $W_L$. Let also 
$\tilde{\sigma} \in {\Spin}^c(W_L)$ be an extension of $\beta(\sigma)\in
{\Spin}^c(V_L)$. Then, the lemma will follow from the fact
that $r_\sigma$ goes to $c(\tilde\sigma)$ by the natural 
map $H^2(W_L,V_L)\rightarrow H^2(W_L)$
provided $\tilde{\sigma}$ is appropriately choosen with respect to $r_\sigma$.
This can be proved undirectly as follows. 
In case when $\sigma$ can be extended to $W_L$,
this is certainly true: indeed, we can take $r_\sigma=0$
and choose as $\tilde\sigma$ the image by $\beta$
of the unique extension of $\sigma$ to $W_L$, so that $c(\tilde\sigma)$ vanishes.
The general case can be reduced to this particular one for 
the following two reasons.
First, it is easily verified that for each Kirby's move $L_1 \to L_2$
between ordered oriented  framed links, the induced bijections 
$\mathcal{S}_{L_1} \rightarrow \mathcal{S}_{L_2}$ and 
$\mathcal{V}_{L_1} \rightarrow \mathcal{V}_{L_2}$, which are
respectively described in \cite[Theorem (I.1)]{Blanchet} and 
Theorem \ref{th:Spin^c_Kirby}, are compatible with the maps 
$\beta:\mathcal{S}_{L_k} \rightarrow \mathcal{V}_{L_k}$ ($k=1,2$)
defined by $\beta\left([r]\right)=\left[B_{L_k}\cdot r\right]$.
Second, according to a theorem of Kaplan \cite{Kap}, 
there exists an oriented framed link $L'$ in
$\mathbf{S}^{3}$ related to $L$ by a finite sequence of Kirby's moves, and
through which $\sigma\in {\Spin}(V_L)$ goes to $\sigma'\in {\Spin}(V_{L'})$
with the property that $\sigma'$ can be extended to $W_{L'}$. 
\end{proof}

\subsubsection{A combinatorial description of $H_2(V_L;\Q/\Z)$}

We maintain the notations used in \S \ref{subsec:discriminant}.
\begin{lemma} 
\label{lem:def_kappa}
There exists a canonical isomorphism
$$\xymatrix{
{\frac{H^{\sharp}}{H}} \ar[r]^-{\kappa}_-{\simeq} & H_2(V_L;\Q/\Z).
}$$
\end{lemma}
\begin{proof}
Consider the following commutative diagram with exact rows and columns:
$$\xymatrix{
&&{0}&{0}\\
{0} \ar[r] & {H_2(V_L;\Q/\Z)} \ar[r] & {H_2(W_L;\Q/\Z)} \ar[r] \ar[u] 
& H_2(W_L,V_L;\Q/\Z) \ar[u]\\
{0} \ar[r]& {H_2(V_L;\Q)} \ar[r] \ar[u]& {H_2(W_L;\Q)} \ar[r]^c \ar[u]^d
& {H_2(W_L,V_L;\Q)} \ar[u]\\
{0} \ar[r] & {H_2(V_L;\Z)} \ar[r] \ar[u] & {H_2(W_L;\Z)} \ar[r] \ar[u]^a 
& {H_2(W_L,V_L;\Z)}\ar[u]^b \\
&&{0}\ar[u]&{0}\ar[u]
}
$$
The group $H^\sharp$ is the subgroup of $H\otimes \Q=H_2(W_L;\Q)$
comprising those $x\in H_2(W_L;\Q)$ 
such that $c(x) \in H_2(W_L,V_L;\Q)$ satisfies
$c(x)\bullet a(y)\in\Z$ for all $y\in H_2(W_L;\Z)$,
where $\bullet$ is the rational intersection pairing in $W_L$. So, we have that
$$H^\sharp=c^{-1}b(H_2(W_L,V_L;\Z)).$$
Seeing $H_2(V_L;\Q/\Z)$ as a subgroup of
$H_2(W_L;\Q/\Z)$, we deduce the announced isomorphism
from the map $d$.
\end{proof}

Recall that the quotient group $H^\sharp/H$, which is denoted by $G_f$ 
in \S \ref{subsec:discriminant}, appears in the short exact sequence (\ref{eq:shortexact}).
We now interpret this sequence as an application of the universal coefficients 
theorem to $V_L$. We denote by $B$ the Bockstein homomorphism associated
to the short exact sequence of coefficients
$$
\xymatrix{
0\ar[r]& \Z \ar[r] &\Q \ar[r] &\Q/\Z \ar[r]& 0.
}
$$
\begin{lemma}
\label{lem:ses}
The  following diagram is commutative:
$$\xymatrix{
{0} \ar[r] & {{\rm{Ker}}\ \widehat{L_f}} \ar[r] \ar[d]_{\kappa|}^\simeq & 
{G_f} \ar[r] \ar[d]_{\kappa}^\simeq& 
{{\rm{Tors}\ \rm{Coker}}\ \widehat{f}} \ar[r] \ar[d]^\simeq_{\varrho|} & 0 \\
{0} \ar[r] & {H_2(V_L)\otimes \Q/\Z} \ar[r] &
{H_2\Big(V_L;\Q/\Z\Big)} \ar[r]_{-P\circ B} & 
{{\rm{Tors}}\ H^2(V_L)} \ar[r]& {0.}
}
$$ 
\end{lemma}

\begin{proof}

It is enough to prove the commutativity of the right square.
Start with a class $m \in H_2\left(V_{L};\Q/\Z\right)$.
It can be written as $m=\left[S\otimes \left[\frac{1}{n}\right]\right]$
where $n$ is a positive integer, $S$ is a $2$--chain in $V_L$ with boundary $\partial S = n \cdot X$
and $X$ is a $1$--cycle. Then, we have that $B(m)=x\in H_1\left(V_L\right)$ if $x$ denotes $[X]$. 
Let also $Y$ be a relative $2$--cycle in $\left(W_L,V_L\right)$ with boundary $\partial Y= X$
and set $y=[Y]\in H_2\left(W_L,V_L\right)$.
Lastly, consider the $2$--cycle $U=n\cdot Y - S$ in $W_L$ and set
$u=[U]\in H=H_2\left(W_L\right)$

Note that $u \otimes \frac{1}{n} \in H\otimes \Q$ 
belongs to the dual lattice $H^\sharp$:
indeed, $P^{-1}\widehat{f}(u)=i_*(u)\in H_2(W_L,V_L)$ equals 
$n\cdot y$ so that $\widehat{f}(u)=n\cdot P(y)$. This also shows that
$\widehat {f_{\Qi}}\left( u \otimes \frac{1}{n} \right)|_H$ $=P(y)$. So,
the map $G_f \rightarrow \hbox{Tors}\ \hbox{Coker}\ \widehat{f}$ that is
featured by the short exact sequence (\ref{eq:shortexact})
sends $\left[ u \otimes \frac{1}{n} \right]$ to $\left[P(y)\right]$.

The canonical map $H\otimes \Q \simeq H_2(W_L;\Q) \rightarrow H_2(W_L;\Q/\Z)$
sends $u \otimes \frac{1}{n}$ to 
$\left[(n\cdot Y-S)\otimes \left[\frac{1}{n}\right]\right] =\left[-S\otimes \left[\frac{1}{n}\right]\right]$.
Consequently, we get that $\kappa\left(\left[ u \otimes \frac{1}{n}\right]\right)=-m$.

The conclusion then follows from the commutativity of the diagram
$$
\xymatrix{
{H_2(W_L,V_L)} \ar[r]^{\partial_*} \ar[d]^P & {H_1(V_L)}\ar[d]^P\\
{H^2(W_L)} \ar[r]^{i^*} & {H^2(V_L),} }
$$
which implies that $\varrho\left([P(y)]\right)=P(x)$.
\end{proof}

\begin{remark}
\label{rem:pairings}
Similarly, the pairing (\ref{eq:evaluation}) can easily be interpreted as the intersection pairing of $V_L$
$$
\xymatrix{
{H_1(V_L)\times H_2(V_L;\Q/\Z)} \ar[r]^-{\bullet} & {\Q/\Z}
}
$$
via the isomorphisms $P^{-1}\varrho: \hbox{Coker}\ \widehat{f} \to H_1(V_L)$
and $\kappa: G_f \to H_2(V_L;\Q/\Z)$.
\end{remark}

\subsection{A $4$--dimensional definition of the linking quadratic function} 

\label{subsec:4def}
Let $M$ be a closed connected oriented $3$--manifold equipped with
a \spinc--structure $\sigma$. In this subsection, we construct the quadratic function
$\phi_{M,\sigma}$ announced in the introduction. 

\begin{lemma} 
\label{lem:4def}
Fix a homology class $m\in H_2(M;\Q/\Z)$.
Consider a quadruplet $(W,\psi, \alpha, w)$ formed by a compact oriented
$4$--manifold $W$, a positive diffeomorphism $\psi: \partial W \to M$, a \spinc--structure
$\alpha$ on $W$ which restricts to $\psi^*(\sigma)$ on the boundary and a class
$w\in H_2(W;\Q)$, the reduction of which in $H_2(W;\Q/\Z)$ coincides with the image of $m$.
Then, the quantity
$$
\phi(M,\sigma,m) = \left[\frac{1}{2}\left(\langle c(\alpha), w \rangle - w\bullet w \right)\right] \in \Q/\Z
$$
does not depend on the choice of such a quadruplet.
\end{lemma}

\begin{remark}
\label{rem:exist}
If  $W$ is a compact oriented $4$--manifold such that $H_1(W)=0$
and there exists a positive diffeomorphism $\psi: \partial W \to M$, 
then the pair $(W,\psi)$ can be completed to a quadruplet $(W,\psi, \alpha, w)$ with the above property.
In particular, such quadruplets do exist since $M$ possesses surgery presentations.
\end{remark}

\begin{proof}
Let $(W',\psi',\alpha',w')$ be another such quadruplet. We wish to compare the 
rational numbers $A:= w \bullet w - \langle c(\alpha), w \rangle$ and
$A':= w'\bullet w' - \langle c(\alpha'), w' \rangle$.

The homology class $m$ of $M$ can be written as $m=\left[S\otimes \left[\frac{1}{n}\right]\right]$,
where $n$ is a positive integer, $S$ is a $2$--chain with boundary $\partial S = n \cdot X$
and $X$ is a $1$--cycle. Then, we have that $B(m)=[X]$. Since the image of $m$ 
in $H_2(W;\Q/\Z)$ belongs to the image of $H_2(W;\Q)$, the image of $[X]\in H_1(M)$ in $H_1(W)$ is zero.
So,  one can find a relative $2$--cycle $Y$ in $(W,\partial W)$ with boundary $\partial Y= \psi^{-1}(X)$.
Consider the $2$--cycle $U=n\cdot Y - \psi^{-1}(S)$ in $W$. Then, by assumption, $w$ can be written
as $w=\left[-U\otimes \frac{1}{n}\right]+w_0 \in H_2(W;\Q)$, where $w_0\in H_2(W;\Q)$ 
belongs to the image of $H_2(W;\Z)$. 
We do the same for $w'$ in $W'$ (getting thus some $Y'$, $U'$, $w_0'$).

Next, we consider the closed oriented $4$--manifold 
$$
\overline{W}:= W \cup_{\psi^{-1} \circ \psi'} (-W').
$$
Gluing rigid \spinc--structures, it is easy to find a \spinc--structure $\overline{\alpha}$ on $\overline{W}$
which restricts to $\alpha$ and $-\alpha'$ on $W$ and $-W'$ respectively.

Set $\overline{Y}= i(Y)-i'(Y')$, where $i$ and $i'$ denote the inclusions of $W$ and $W'$ respectively.
This is a $2$--cycle in $\overline{W}$ with the property that the identity
$$
\left[\overline{Y} \otimes 1\right]
= \left[i(U)\otimes 1/n - i'(U') \otimes 1/n \right]\\
=\left( -i_*(w) + i_*(w_0) \right) + \left( i'_*(w') -i'_*(w'_0)\right)
$$
holds in $H_2\left(\overline{W};\Q\right)$.
It follows from this identity that
\begin{equation}
\label{eq:double_intersection}
\begin{array}{rcl}
[\overline{Y}] \bullet [\overline{Y}] &=& 
\left(w \bullet w + w_0 \bullet w_0 - 2\cdot w \bullet w_0\right)\\
&&+ \left(-w' \bullet w' - w_0' \bullet w_0' + 2\cdot w' \bullet w_0' \right)
\end{array}
\end{equation}
and that
\begin{equation}
\label{eq:evaluations}
\left\langle c\left(\overline{\alpha}\right), [ \overline{Y}] \right\rangle
= \left(-  \langle c(\alpha), w \rangle + \langle c(\alpha), w_0 \rangle \right)+
\left( \langle c(\alpha'), w' \rangle - \langle c(\alpha'), w_0' \rangle\right).
\end{equation}

Recall that $w_0 \in H_2(W;\Q)$ and $w'_0 \in H_2(W';\Q)$ come from integral classes.
Then, by the Wu formula and the fact that a Chern class reduces modulo $2$ 
to the second Stiefel--Whitney class,
the integers $[\overline{Y}] \bullet [\overline{Y}]$,
$w_0\bullet w_0$ and $w_0'\bullet w_0'$ are congruent modulo $2$ to
$\left\langle c\left(\overline{\alpha}\right), [\overline{Y}] \right\rangle$,
$\langle c(\alpha), w_0 \rangle$ and $\langle c(\alpha'), w_0' \rangle$ respectively.
Adding (\ref{eq:double_intersection}) to (\ref{eq:evaluations}), we find that
$$
A - A'  - 2\cdot w \bullet w_0 + 2\cdot w' \bullet w_0' \equiv 0 \ \ \hbox{mod } 2.
$$
Because the image of $w\in H_2(W;\Q)$ in $H_2(W;\Q/\Z)$ comes from
$H_2(M;\Q/\Z)$ and because $w_0\in H_2(W;\Q)$ comes from $H_2(W;\Z)$, the rational
number  $w \bullet w_0$ belongs to $\Z$. The same holds for $w' \bullet w_0'$.
We conclude that the rational number $A-A'$ belongs to $2\cdot \Z$.
\end{proof}

\begin{remark}
A universal class $u\in H^1\left(K(\Q/\Z,1);\Q/\Z\right)$ induces a homomorphism
$$
\xymatrix{
\Omega_3^{{\rm \scriptsize Spin}^c}\left(K(\Q/\Z,1)\right) \ar[r] & \Q/\Z}
$$
defined by $\left[(M,\sigma,f)\right] \mapsto \phi\left(M,\sigma,P^{-1}f^*(u)\right)$.
This follows from the definition of $\phi$ in Lemma \ref{lem:4def}.
\end{remark}

Consider the linking pairing $\lambda_M: \hbox{Tors }H_1(M) \times \hbox{Tors }H_1(M) \rightarrow \Q/\Z$. 
Composing this with the Bockstein $B$, one gets a symmetric bilinear pairing 
$$
\xymatrix{H_2(M;\Q/\Z) \times H_2(M;\Q/\Z) \ar[r]^-{L_M}& \Q/\Z}
$$ 
with radical $H_2(M)\otimes\Q/\Z$. 
Using a cobordism $W$ as in Remark \ref{rem:exist}, one easily proves,
for any $m,m' \in H_2(M;\Q/\Z)$, the following identity:
$$
\phi(M,\sigma,m+m') - \phi(M,\sigma,m) - \phi(M,\sigma,m')= m \bullet B(m') = L_M(m,m'). 
$$
\begin{definition}
\label{def:quad_spinc}
The \emph{linking quadratic function} of the \spinc--manifold $(M,\sigma)$ is the map denoted by
$$
\xymatrix{ H_{2}\left(M;\Q/\Z\right) \ar[r]^-{\phi_{M,\sigma}} & \Q/\Z}
$$
and defined by $m\mapsto \phi(M,\sigma,m)$. 
\end{definition}

The discriminant construction allows us to compute combinatorially
the quadratic function $\phi_{M,\sigma}$, as soon as
a surgery presentation of the \spinc--manifold $(M,\sigma)$ is given. Indeed,
let $L$ be an ordered oriented framed link in $\mathbf{S}^{3}$
together with a positive diffeomorphism $\psi: V_L \rightarrow M$.
With the notations from \S \ref{subsec:combinatorial}, $(H,f)$ still denotes the bilinear lattice 
$\left(H_2(W_L), \textrm{intersection pairing of }W_L\right)$, to which 
the constructions from \S \ref{subsec:discriminant} apply.
Let also $c\in \hbox{Char}(f)$ represent $\psi^*(\sigma)\in \hbox{Spin}^c(V_L)$
(in the sense of Lemma \ref{lem:combinatorial_spinc}).
Then, as can be verified from the definitions, the following diagram commutes:
\begin{equation}
\label{eq:triangle}
\xymatrix{
H_2(M;\Q/\Z) \ar[rr]^-{\phi_{M,\sigma}}  && \Q/\Z\\
H_2(V_L;\Q/\Z) \ar[u]_-\simeq^-{\psi_*}&&\\
G_f \ar[u]^-{\kappa}_-\simeq \ar[rruu]_-{-\phi_{f,c}}&&
}
\end{equation}
Note that, in this context, the pairings $\lambda_f$ and $L_f$ are 
topologically interpreted as $-\lambda_M$ and $-L_M$ respectively.

\subsection{Properties of the linking quadratic function}

In this subsection, we fix a closed connected oriented  $3$--manifold $M$ and
prove properties of the map $\phi_M: \Spinc(M) \to \Quad(L_M)$ defined by $\sigma \mapsto \phi_{M,\sigma}$.
Those properties are proved ``combinatorially''  using (\ref{eq:triangle}), 
but may also be proved directly from the very definition of $\phi_{M,\sigma}$.

Next lemma says that $\phi_{M,\sigma}$ is determined
on $H_2(M)\otimes \Q/\Z$ by the Chern class $c(\sigma)$. 
Recall that the modulo $2$ reduction of $c(\sigma)$ is $w_2(M)=0$.

\begin{lemma}
\label{lem:behaviour_on_the_kernel}
For any $\sigma \in {\rm{Spin}}^c(M)$, the function $\phi_{M,\sigma}$ is linear on $H_2(M)\otimes \Q/\Z$:
\begin{displaymath}
\forall x\otimes [r] \in H_2(M)\otimes \Q/\Z,\ \ 
\phi_{M,\sigma}\left(x\otimes [r]\right)
=\frac{\langle c(\sigma),x \rangle}{2} \cdot [r] \in \Q/\Z.
\end{displaymath} 
\end{lemma}
\begin{proof}
The first statement follows from the fact that 
$\hbox{Ker}\ \widehat{L_M}=H_2(M)\otimes \Q/\Z$. As for the second statement,
it suffices to prove it when $M=V_L$. Suppose that 
$\sigma$ is represented by the characteristic form $c\in \hbox{Char}(f)$ 
and that $x\in H_2(V_L)$ goes
to $y$ in $H=H_2(W_L)$. Then, $x\otimes [r]$
as an element of $H_2(V_L; \Q/\Z)$ corresponds to $[y \otimes r]$
in $H^\sharp/H$. Consequently, we have that $\phi_{M,\sigma}\left(x\otimes [r]\right)=$
$-\phi_{f,c}\left([y \otimes r]\right)$
$=-\frac{1}{2}\left(r^2f(y,y)-r\cdot c(y)\right)$ mod $1$. Since $y$ belongs to $\hbox{Ker}\ \widehat{f}$,
we obtain that $\phi_{M,\sigma}\left(x\otimes [r]\right)=$
$\frac{1}{2} r \cdot c(y) \hbox{ mod } 1= \frac{1}{2} r\cdot \langle c(\sigma),x \rangle  \hbox{ mod } 1$,
by Remark \ref{rem:combinatorial_Chern}.
\end{proof}

Let us consider, for a while, the case when $\sigma \in {\Spin}^c(M)$ is torsion.
Then, Lemma \ref{lem:behaviour_on_the_kernel} implies that
$\phi_{M,\sigma}$ vanishes on $H_2(M)\otimes\Q/\Z$: Consequently, 
$\phi_{M,\sigma}$ factors to a quadratic function over $\lambda_M$.
In this torsion case, our linking quadratic function is readily seen to 
agree with that of \cite{De} and, up to a minus sign, with that of \cite{Gille}.
In the next subsection, it is also shown to coincide with that of \cite{LW}.

In particular, $\sigma$ may arise from a \spin--structure on $M$,
which happens if and only if $c(\sigma)$ vanishes.
Then, the factorization of $\phi_{M,\sigma}$ to $\hbox{Tors}\ H_1(M)$ coincides
with the linking quadratic form defined in \cite{LL}, \cite{MS} or
\cite{TCohomology}. In \cite{Mas}, this quadratic form is used
to classify degree $0$ invariants in the \spin--refinement 
of the Goussarov--Habiro theory.\\

In the sequel, we will use the homomorphism
$$
\xymatrix{
{H^2(M)} \ar[r]^-{\mu_M} & 
{\hbox{Hom}\left(H_2(M;\Q/\Z), \Q/\Z\right)}
}
$$
defined by $\mu_M(y)=\left\langle y,-\right\rangle$. 

\begin{lemma} \label{lem:homogeneity-defect}
For any $\sigma\in {\rm{Spin}}^c(M)$, the Chern class $c(\sigma)$ is sent by $\mu_M$ 
to the homogeneity defect $d_{\phi_{M,\sigma}}: H_2(M;\Q/\Z) \to \Q/\Z$
of the quadratic function $\phi_{M,\sigma}$.
\end{lemma}
\begin{proof}
Again suppose that $M=V_L$ 
and that $\sigma$ is represented by $c\in \hbox{Char}(f)$. 
Take $x\in H_2(V_L;\Q/\Z)$ represented
by $y \in H^\sharp$. One computes that $\phi_{M,\sigma}(x)-\phi_{M,\sigma}(-x)$
$=-\phi_{f,c}\left([y]\right)+\phi_{f,c}\left(-[y]\right)$
$=c_{\Qi}(y) \hbox{ mod } 1$
$=\langle c(\sigma),x \rangle$, by Remark \ref{rem:combinatorial_Chern}.
\end{proof}
Recall that $\hbox{Spin}^c(M)$ is an affine space over $H^2(M)$
and that $\Quad(L_M)$ is an affine space over 
$\hbox{Hom}\left(H_2(M;\Q/\Z), \Q/\Z\right)$.  Let
\begin{displaymath}
\xymatrix{
{\hbox{Hom}\left(H_2(M),\Z\right)} \ar[r]^{j_M\quad \quad} &
{\hbox{Hom}\left(H_2(M)\otimes \Q/\Z,\Q/\Z\right)}
}
\end{displaymath}
be the homomorphism defined by $j_M(l)=l\otimes \Q/\Z$.
Next result, which contains Theorem \ref{th:embedding},
is a direct application of Theorem \ref{th:algebraic-embedding}
and Remark \ref{rem:pairings}.
\begin{theorem}
\label{th:phi_M}
The map $\phi_M:{\rm{Spin}}^c(M)\rightarrow {\rm{Quad}}(L_M)$ 
is an affine embedding over the group monomorphism $\mu_M$.
Moreover, a function $q\in {\rm{Quad}}(L_M)$ belongs to ${\rm{Im}}\ \phi_{M}$ if and only if
$q|_{H_2(M)\otimes \Qi/\Zi}$ belongs to ${\rm{Im}}\ j_M$.
\end{theorem}
\begin{remark}
The map $\phi_M$ is bijective if and only if $M$ is a rational homology $3$--sphere.
\end{remark}

\subsection{An intrinsic definition of the linking quadratic function} 

\label{subsec:geometric} 
Let $M$ be a closed connected oriented $3$--manifold equipped with
a \spinc--structure $\sigma$. In this subsection, we give for
the quadratic function $\phi_{M,\sigma}$
an intrinsic formula  which does not refer to $4$--dimensional cobordisms.

Here is the idea. Take a $x\in H_2(M;\Q/\Z)$.
It follows from Lemma \ref{lem:homogeneity-defect} that 
$$ 2\cdot \phi_{M,\sigma}(x)= 
L_M(x,x)+ \langle c(\sigma),x \rangle \in \Q/\Z.$$ 
For any $y\in\Q/\Z$, we denote by $\frac{1}{2}\cdot y$
the set of elements $z$ of $\Q/\Z$ such that $z+z=y$. 
We are going to select, \emph{correlatively}, an element 
$z_1$ in $\frac{1}{2}\cdot L_M(x,x)$ \emph{and} an element  $z_2$ in 
$\frac{1}{2}\cdot \langle c(\sigma),x \rangle$ such that $\phi_{M,\sigma}(x)=z_1+z_2$.

Write $x\in H_2(M;\Q/\Z)$ as $x=\left[ S\otimes \left[1/n\right] \right]$,
where $n$ is a positive integer and $S$ is an oriented immersed surface 
in $M$ with boundary $n\cdot K$,
a bunch of $n$ parallel copies of an oriented knot $K$ in $M$.
Apply now the following stepwise procedure:
\begin{itemize}
\item[$\centerdot$] \emph{Step 1.} 
Choose a nonsingular vector field $v$ on $M$ representing
$\sigma$ as an Euler structure, 
and which is transverse to $K$ (we claim that it is possible to find such $v$).
\item[$\centerdot$] \emph{Step 2.} Let $V$ be a sufficiently small regular neighborhood
of $K$ in $M$ and let $K_v$ be the parallel of $K$, lying on $\partial V$, 
obtained by pushing $K$ along the trajectories of $v$. By an isotopy, 
ensure that $S$ is in transverse position with respect to $K_v$ 
with boundary contained in the interior of $V$.
\item[$\centerdot$] \emph{Step 3.}
Define a \spin--structure $\alpha_v$ on $\partial\left(M\setminus
\hbox{int}\left(V\right)\right)$
by requiring its Atiyah--Johnson quadratic form $q_{\alpha_v}$ 
(\S \ref{subsubsec:relative_Chern}) to be such that
$$q_{\alpha_v}\left(\left[\textrm{meridian of } K\right]\right)=0 \ \textrm{ and } \
\ q_{\alpha_v}\left(\left[K_v\right]\right)=1.$$ 
\item[$\centerdot$] \emph{Step 4.} Together with the vector field tangent to $K_v$, 
$v$ represents a \spinc--structure $\sigma_v$ on $M\setminus \hbox{int}\left(V\right)$ relative 
to the \spin--structure $\alpha_v$ (we claim this). Consider the Chern class
$c(\sigma_v) \in H^2\left(M\setminus \hbox{int}\left(V\right),
\partial\left(M\setminus \hbox{int}\left(V\right)\right)\right).$
\end{itemize}
\begin{proposition}
\label{prop:intrinsic_formula}
By applying the above procedure, we get
\begin{equation}
\label{eq:intrinsic_formula}
\phi_{M,\sigma}(x) = \underbrace{\left[\frac{1}{2n}\cdot K_v\bullet S\right]}
_{\in\frac{1}{2}\cdot L_M(x,x)} + \underbrace{\left[\frac{1}{2n}\cdot 
\left\langle c(\sigma_v),\left[S\cap \left(M\setminus{\rm{int}}\left(V
\right)\right)\right]\right\rangle +\frac{1}{2}\right]}_{\in
\frac{1}{2}\cdot \langle c(\sigma),x\rangle} 
\in \Q/\Z.
\end{equation} 
\end{proposition}

In \cite{LW}, Looijenga and Wahl associate a quadratic function over $\lambda_M$
to each pair $(M,\mathcal{J})$ formed by 
\begin{itemize}
\item[$\centerdot$] a closed connected oriented $3$--manifold $M$,
\item[$\centerdot$] a homotopy class of complex structures $\mathcal{J}$
on $\R \oplus \hbox{T} M$ whose first Chern class is torsion.
\end{itemize}
There is a \spinc--structure $\omega\left(\mathcal{J}\right)$ associated to $\mathcal{J}$
(see \S \ref{subsubsec:U}). By assumption, its Chern class is torsion so that
$\phi_{M,\omega\left(\mathcal{J}\right)}$ factors to a quadratic function over 
$\lambda_M$. One can verify, using the inverse of $\omega$ described in the proof
of Lemma \ref{lem:U_dim3}, that formula (\ref{eq:intrinsic_formula}) is equivalent 
in this case to formula (3.4.1) in \cite{LW}. 

\begin{proof}[Proof of Proposition \ref{prop:intrinsic_formula}]
First of all, we have to justify that the above procedure can actually
be carried out.

We begin by proving the claim of Step 1. Let $v$ be an arbitrary nonsingular
vector field on $M$ representing $\sigma$. It suffices to prove
the following claim.
\begin{claim}
Let $w$ be an arbitrary nonsingular vector field tangent to $M$ defined on $K$.
Then, $v$ can be homotoped so as to coincide with $w$ on $K$.
\end{claim}
\begin{proof}
Choose a tubular neighborhood $W$ of $K$, plus
an identification $W=(2\mathbf{D}^2)\times \mathbf{S}^1$ such that
$K$ corresponds to $0\times \mathbf{S}^1$. We denote
by $\left(e_1,e_2\right)$ the standard basis of $\R^2\supset 2\mathbf{D}^2$.
We define $\pi:W\rightarrow K$ to be the projection on the core.
The solid torus $W$ is parametrized by the cylindric coordinates 
\begin{displaymath}
\left(\left(r\in [0,2],\theta\in \R/ 2\pi\Z \right), \phi \in \R/2\pi\Z\right).
\end{displaymath}
If $p,q \in W$ are such that $\pi(p)=\pi(q)$ 
(i.e., they belong to the same meridional disk $2\mathbf{D}^2\times *$), 
we define the transport map $t_{p,q}:\hbox{T}_p W \to \hbox{T}_q W$ 
as the unique linear map fixing
the basis $\left(e_1,e_2,\frac{\partial}{\partial \phi}\right)$. 
Deform the vector field $v$ through the homotopy 
$\left(v^{(t)}\right)_{t\in[0,1]}$ given at time $t$ and point $p\in W$ by
\begin{displaymath} 
v^{(t)}_{p}= \left\{ 
\begin{array}{ll}
t_{\pi(p),p}\left(v_{\pi(p)}\right) & \textrm{ if } r(p) \in [0,t]\\
t_{q(p,t),p}\left(v_{q(p,t)}\right) & \textrm{ if } r(p)\in [t,2],
\textrm{ with } q(p,t)=\left(\frac{r(p)-t}{1-t/2},\theta(p),\phi(p)\right)
\end{array} \right.
\end{displaymath}
and at time $t$ and point $p\notin W$ by $ v^{(t)}_{p}=v_p$. 
After such a deformation, the vector field $v$ 
satisfies the following property:
$\forall p\in \mathbf{D}^2\times \mathbf{S}^1, 
\ t_{p,\pi(p)}(v_p)=v_{\pi(p)}$. Now, since 
$\pi_1(\mathbf{S}^2)$ is trivial, 
$v|_K$ and $w$ have to be homotopic;
let $\left(w^{(t)}\right)_{t\in[0,1]}$ be such a homotopy, 
beginning at $w^{(0)}=v|_K$ and ending at $w^{(1)}=w$. 
The homotopy given by
\begin{displaymath}
v^{(t)}_{p}= \left\{
\begin{array}{lr}
t_{\pi(p),p}\left(w^{(t-r(p))}_{\pi(p)}\right) & \textrm{ if } r(p) \in [0,t]\\
v_p & \textrm{ if } r(p) \in [t,2]
\end{array} \right.          
\end{displaymath}
if $p \in W$ and by $v^{(t)}_{p}=v_p$ if $p\notin W$, allows us to deform
$v$ to a nonsingular vector field which coincides with $w$ on $K$. 
\end{proof}
\noindent
Since $v$ is now transverse to $K$, we can find a regular 
neighborhood $V$ of $K$ in $M$ plus an identification
$V=\mathbf{D}^2\times \mathbf{S}^1$, such that $K$ corresponds
to $0\times \mathbf{S}^1$ and such that $v|_V$ corresponds to $e_1$
(recall that $(e_1,e_2)$ denotes the standard basis of 
$\R^2\supset \mathbf{D}^2$). We apply steps 2 and 3 (note that $K_v$
then corresponds to $1\times \mathbf{S}^1$) and 
we now prove the claim of Step 4. 
Let $\tau_v\in \hbox{Spin}(V)$ be defined by the trivialization 
$\left(e_1,e_2,\frac{\partial}{\partial \phi}\right)$ of $\hbox{T}V$.
Since $\left(\tau_v|_{ \partial V}\right)|_{1\times \mathbf{S}^1}$ 
is the non-bounding \spin--structure and since 
$\left(\tau_v|_{ \partial V}\right)|_{\partial \mathbf{D}^2 \times 1}$
spin bounds, we have that $\tau_v|_{ \partial V}=-\alpha_v$,
i.e. $\tau_v$ belongs to $\Spin(V,-\alpha_v)$ with the notation
of Remark \ref{rem:from_spin_to_spinc_relative}. Thus, $v|_{ M\setminus \hbox{\scriptsize int}\left(V\right) }$ 
together with the trivialization $\left.\left(e_1,e_2,\frac{\partial}{\partial \phi}\right)\right|_{\partial V}$ 
of $\left.\hbox{T} \left( M\setminus \hbox{int}\left(V\right) \right)\right|_{\partial V}$ 
define a $\sigma_v\in \Spinc\left(M\setminus \hbox{int}\left(V\right),\alpha_v\right)$,
as claimed in Step 4. For further use, note that $\sigma$ 
is the gluing $\sigma_v\cup\beta(\tau_v)$, where 
$\beta:\Spin(V,-\alpha_v) \to \Spin^c(V,-\alpha_v)$
has been defined in Remark \ref{rem:from_spin_to_spinc_relative}.\\
 
Set $z_1=\left[1/2n \cdot K_v\bullet S\right] \in \Q/\Z$
and $z_2=\left[1/2n \cdot \left\langle c(\sigma_v),[S']\right\rangle 
+1/2\right]\in \Q/\Z$, where $S'=S\cap \left(M\setminus
\hbox{int}\left(V\right)\right)$. We  have that
$$
2\cdot z_1=\left[1/n \cdot K_v\bullet S\right]=
\left[\lambda_M\left(B(x),B(x)\right)\right]=L_M(x,x).
$$
Moreover, we have that
\begin{eqnarray*} 2\cdot z_2&=&
\left[1/n \cdot \left\langle c(\sigma_v),[S']\right\rangle\right]\\
&=&\left[1/n \cdot P^{-1}\left(c(\sigma_v)\right)\bullet [S']\right]
\ (\textrm{intersection in } M\setminus \hbox{int}(V))\\
&=& P^{-1}\left(c(\sigma)\right) \bullet x 
\ (\textrm{intersection in } M)\\
&=& \langle c(\sigma), x \rangle
\end{eqnarray*}
where the third equality follows from the facts that 
$x=\left[S\otimes \left[1/n\right]\right]$,
$P^{-1}\left(c(\sigma)\right)=i_* P^{-1}\left(c(\sigma_v)\right)+
i_* P^{-1}\left(c(\beta(\tau_v))\right)\in H_1(M)$ 
(since $\sigma=\sigma_v\cup\beta(\tau_v)$) and 
$c(\beta(\tau_v))=0$ (by Remark \ref{rem:relative_Chern}).\\

We now prove formula (\ref{eq:intrinsic_formula}), i.e. , the equality
$\phi_{M,\sigma}(x)=z_1+z_2$. Let us work with surgery presentations
(even if we could use more general cobordisms as well).
Let $M'$ be the $3$--manifold obtained from $M$ by doing surgery along the framed
knot $\left(K,(e_1,e_2)\right)$. Conversely, $M$ is the result of the surgery 
on $M'$ along the dual knot $K'$ of $K$.
Pick a surgery presentation $V_{L'}$ of $M'$; up
to isotopy, the knot $K'\subset M'$ is in $\mathbf{S}^3\setminus L'$.
We then find a surgery presentation 
$V_L$ of $M$ by setting $L$ to be $L'$ union $K'$ 
with the appropriate framing. This surgery presentation of $M$ has
the following advantage: $K$ bounds in the trace $W_L$ of the surgery
a disk $D$ whose normal bundle is trivialized by some extension of the
trivialization $(e_1,e_2)$ of the normal bundle of $K$ in $M$.
We use the notations fixed in \S \ref{subsec:combinatorial}. In particular,
$H=H_2(W_L)$ and $f:H\times H\to \Z$ 
is the intersection pairing of $W_L$. We define the $2$--cycle
$U=n\cdot D-S$ where $n\cdot D$ is a bunch of $n$ parallel copies of the
disk $D$ with boundary $n\cdot K$; we also set $u=[U] \in H$. 
Then $u\otimes \frac{1}{n}$ belongs
to $H^\sharp$ and the isomorphism $\kappa: H^\sharp/H \rightarrow H_2(M;\Q/\Z)$ 
sends $\left[u \otimes \frac{1}{n}\right]$ to $-x=-\left[S\otimes \left[\frac{1}{n}\right]\right]$ 
(see the proof of Lemma \ref{lem:ses}).
So, by diagram (\ref{eq:triangle}), we obtain that
$$\phi_{M,\sigma}(x)=-
\phi_{f,c}\left(-\left[  u \otimes \frac{1}{n} \right]\right)=
-\frac{1}{2}\left( \frac{1}{n^2}f(u,u)+\frac{1}{n}c(u)\right) \hbox{ mod } 1$$
where $c$ is a characteristic form representative for $\sigma$.

We calculate the quantity $f(u,u)$. The $2$--cycle
$U$ is a representant of $u$. Let $D'$ be a push-off of $D$ 
by the extension of $e_1=v|_V$ in such a way that $\partial D'$ is $K_v$.
Let also $A$ be the annulus of an isotopy from
$-K_v$ to $K$ in $V$ (e.g. $A=-[0,1]\times \mathbf{S}^1$ in 
$V=\mathbf{D}^2\times \mathbf{S}^1$). A second representative for $u$
is $U'=n\cdot D'+n\cdot A-S$. 
By adding a collar to $W_L$ and stretching the top of $U'$, we can make $U$ in 
transverse position with $U'$
(see Figure \ref{fig:transverse}).
\begin{figure}[!h]
\begin{center}
\includegraphics[width=10cm,height=6cm]{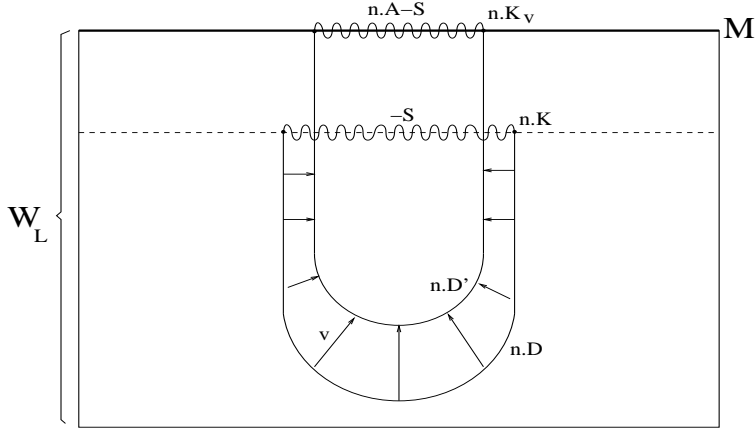}
\caption{Two representants of $u$ in transverse position.}
\label{fig:transverse}
\end{center}
\end{figure}
So, we have that $f(u,u)= U\bullet U'=-n S\bullet K_v$ where the first intersection
is calculated in $W_L$ and the second one in $M$; we are led to
\begin{equation}
\label{eq:intrinsic_formula_intermediate}
\phi_{M,\sigma}(x)=\frac{1}{2n}S\bullet K_v -\frac{1}{2n}c(u) \hbox{ mod } 1.
\end{equation}

We are now interested in the quantity $c(u)$. 
Let $\tilde\sigma$ be an extension of $\sigma$ to the manifold $W_L$ and
let $\xi$ be the isomorphism class of principal $\U(1)$--bundles 
on $W_L$ defined by $\tilde\sigma$; 
then $c$ can be choosen to be $c_1(\xi)$. 
Let $p$ be a representant of $\xi$ and let tr be a trivialization
of $p$ on $\partial V$. 
Decompose the singular surface $U'$ as $U'=U'_1\cup U'_2 \cup U'_3$,
where $U'_1=n\cdot D'$, $U'_2=n\cdot A \cup \left(-S\cap V\right)$
and $U'_3=-S'$. By desingularizing $U'$ so as to
be reduced to a calculus of obstructions in an oriented manifold, we obtain that
\begin{equation}
\label{eq:calcul_c(u)}
c(u)=\left\langle c_1\left(p|_{U'}\right),[U']\right\rangle=
\sum_{i=1}^3 \left\langle c_1\left(p|_{U_i'},
\hbox{tr}|_{\partial U'_i}\right),[U'_i]\right\rangle\in \Z,
\end{equation}
where $c_1\left(p|_{U_i'},\hbox{tr}|_{\partial U'_i}\right) \in
H^2\left(U'_i,\partial U'_i\right)$ is the obstruction to extend
the trivialization $\hbox{tr}|_{\partial U'_i}$ of $p|_{U_i'}$
on $\partial U'_i$ to the whole of $U'_i$. 
Let $V'\subset W_L$ be the solid torus such that 
$M'=M\setminus \hbox{int}(V)\cup V'$. 
For an appropriate choice of $\tilde{\sigma}$,
there exists a \spinc--structure $\sigma_1\in \hbox{Spin}^c(V',-\alpha_v)$
such that $\sigma_v\cup \sigma_1={\tilde{\sigma}}|_{M'}$.
Also, for some appropriate choices of $p$ in the class $\xi$ and tr, we have
\begin{eqnarray*}
c_1\left(p|_{V'},\hbox{tr}\right)&=&c(\sigma_1) \in H^2(V',\partial V'),\\
c_1\left(p|_{V},\hbox{tr}\right)&=&c(\beta(\tau_v)) \in H^2(V,\partial V),\\
c_1\left(p|_{M\setminus \hbox{int}(V)},\hbox{tr}\right)&=&
c(\sigma_v) \in H^2\left(M\setminus \hbox{int}(V),\partial 
\left(M\setminus \hbox{int}(V) \right)\right).
\end{eqnarray*}
Then, equation (\ref{eq:calcul_c(u)}) becomes
$$c(u)=n\cdot \langle c(\sigma_1),[D']\rangle+
\langle c(\beta(\tau_v)), [U'_2]\rangle-
\langle c(\sigma_v), [S'] \rangle \in \Z.$$
From the fact that $c(\beta(\tau_v))=0$, we deduce that
$$\frac{1}{2n}\cdot c(u)= -\frac{1}{2n}\cdot
\langle c(\sigma_v), [S'] \rangle +
\frac{1}{2}\cdot \langle c(\sigma_1),[D']\rangle\in \Q.$$
Then, showing that $\langle c(\sigma_1),[D']\rangle$ is an odd
integer together with (\ref{eq:intrinsic_formula_intermediate})
will end the proof of the proposition. Since 
$\langle c(\sigma_1),[D']\rangle=
q_{-\alpha_v}\left(\partial_* [D']\right)=
q_{\alpha_v}\left(\left[K_v\right]\right)=1$ mod 2 (by Lemma \ref{lem:Chern_mod2}),
we are done.
\end{proof}

\section{Goussarov--Habiro theory for three--manifolds with complex spin structure} 

\label{sec:FTI}
In this section, we explain how the Goussarov--Habiro theory 
can be extended to the context of $3$--manifolds equipped with a \spinc--structure.
Then, using the linking quadratic function, 
we prove Theorem \ref{th:Matveev_spinc} stated in the introduction.
This amounts to identifying the degree $0$ invariants in the generalized theory.

\subsection{Review of the $Y$--equivalence relation}

\label{subsec:review}
Recall that the Goussarov--Habiro theory is a theory of finite type invariants
for compact oriented $3$--manifolds \cite{Goussarov,Habiro,GGP}
and is based on the $Y$--surgery as elementary move . 
In this subsection, we just recall how this surgery move is defined.

\begin{figure}[h!]
\begin{center}
\includegraphics[height=5cm,width=5cm]{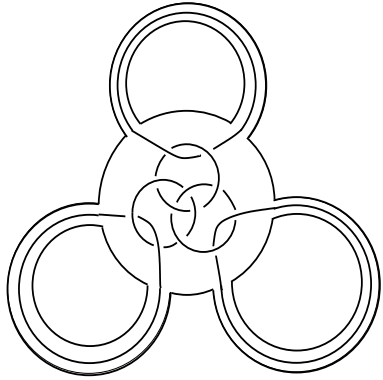} 
\caption{The framed link $B$.}
\label{fig:B}
\end{center}
\end{figure}

Suppose that $M$ is a compact oriented $3$--manifold.
Let $j: H_3 \hookrightarrow M$ be a positive
embedding of the genus $3$ handlebody into the interior of $M$. Set
\begin{displaymath}
M_{j}= M \setminus \hbox{int}\left(\hbox{Im}(j)\right)
\cup_{j|_{\partial H_3}} \left(H_{3}\right)_{B}.
\end{displaymath}
Here, $\left(H_{3}\right)_{B}$ is the surgered handlebody
along the six--component framed
link $B$ shown on Figure \ref{fig:B} with the blackboard framing convention.
\begin{remark}
\label{rem:Y_as_a_twist}
Observe that there is a canonical inclusion 
$M\setminus \hbox{int}\left(\hbox{Im}(j)\right) \hookrightarrow M_j$.
One can define a self-diffeomorphism $h$ of $\partial H_3$
(explicitely, as the composition of $6$ Dehn twists) such that there exists a diffeomorphism
\begin{equation}
\label{eq:gluing}
M_j \cong M\setminus \hbox{int}\left(\hbox{Im}(j)\right)
\cup_{j|_{\partial H_3}\circ h} H_3
\end{equation}
restricting to the identity on $M\setminus \hbox{int}\left(\hbox{Im}(j)\right)$.
Moreover, $h$ can be verified to act trivially in homology.
\end{remark}

\begin{figure}
\begin{center}
\includegraphics[height=3.5cm,width=4cm]{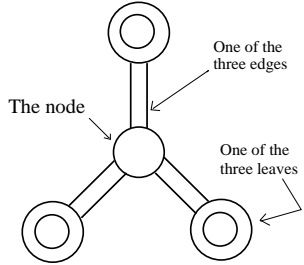} 
\caption{A $Y$--graph.}
\label{fig:Y}
\end{center}
\end{figure}

A \emph{$Y$--graph} $G$ in $M$ is an embedding of the surface drawn in Figure \ref{fig:Y}
into the interior of $M$. This surface, of genus $0$ with
$4$ boundary components, is decomposed
between \emph{leaves}, \emph{edges} and \emph{node}.
Let $j:H_3 \hookrightarrow M$ be a trivialization 
of a regular neighborhood of $G$ in $M$. 
The embedding $j$ is unique, up to ambiant isotopy.
\begin{definition}
The manifold obtained from $M$ by \emph{$Y$--surgery along $G$},
denoted by $M_G$, is the positive diffeomorphism class of the manifold $M_j$.
The \emph{$Y$--equivalence} is the equivalence relation 
among compact oriented $3$--manifolds generated by $Y$--surgeries
and positive diffeomorphisms.
\end{definition}
\begin{remark}
The $Y$--surgery move has been introduced by Goussarov \cite{Goussarov}
and is equivalent to Habiro's ``$A_1$--move'' \cite{Habiro}. 
It is equivalent to Matveev's ``Borromean surgery''  as well,
hence the $Y$--equivalence relation is characterized in \cite{Matveev}. 
\end{remark} 

\subsection{The $Y^c$--equivalence relation}

We define the $Y^c$--surgery move announced in the introduction,
and we outline how this suffices to extend the Goussarov--Habiro theory 
to manifolds equipped with a \spinc--structure.

\subsubsection{Twist and \spinc--structures}

\label{subsubsec:twist}

As in \S \ref{subsec:gluing}, we consider a closed oriented $3$--manifold
\begin{displaymath}
M=M_1 \ \cup_f \ M_2 
\end{displaymath}
obtained by gluing two compact oriented $3$--manifolds $M_1$ and $M_2$
with a positive diffeomorphism $f:-\partial M_2 \to \partial M_1$. We add the
assumption that \emph{$\partial M_2$ is connected.}

Let $h: \partial M_2 \to \partial M_2$ be a diffeomorphism
\emph{which acts trivially in homology} and consider the manifold 
$$
M'=M_1 \ \cup_{f\circ h} \ M_2.
$$
The manifold $M'$ is said to be obtained from $M$ by a \emph{twist}.
By Remark \ref{rem:Y_as_a_twist},
the $Y$--surgery move is an instance of a twist move. 

By a Mayer--Vietoris argument, there is an isomorphism $\Phi:H_1(M)\to H_1(M')$
which is unambiguously defined by the commutative diagram
\begin{displaymath}
\xymatrix{
& H_1(M)\ar[dd]^{\Phi}_\simeq &\\
H_1(M_1) \ar[ru]^{j_{1,*}}\ar[rd]_{j'_{1,*}} & &
H_1(M_2)\ar[lu]_{j_{2,*}}\ar[ld]^{j'_{2,*}}\\
& H_1(M') &}
\end{displaymath}
where $j_1$, $j_2$, $j'_1$ and $j'_2$ denote inclusions.
\begin{proposition} 
\label{prop:twist}
The twist from $M$ to $M'$ induces a canonical bijection
\begin{displaymath}
\xymatrix{{\rm{Spin}}^c(M) \ar[r]^-{\Omega}_-{\simeq} & {\rm{Spin}}^c(M')}
\end{displaymath}
which is affine over $P \Phi P^{-1}: H^2(M) \to H^2(M')$. Moreover, the diagram 
\begin{displaymath}
\xymatrix{
{\rm{Spin}}^c(M)\ar[r]^-\Omega\ar[d]_c & {\rm{Spin}}^c(M')\ar[d]^c \\
H^2(M) \ar[r]_-{P \Phi P^{-1}} & H^2(M')
}
\end{displaymath}
is commutative.
\end{proposition}
\begin{proof}
For any $\alpha \in \Spin^c(M)$, we define $\Omega(\alpha)$ as follows.
\emph{Choose} $\sigma_2 \in \hbox{Spin}\left(\partial M_2\right)$ and set
$\sigma_1=f_*(-\sigma_2)\in \Spin(\partial M_1)$. Since 
$h_*:H_1(\partial M_2;\Z_2)\to H_1(\partial M_2;\Z_2)$ is the identity,
$h$ acts trivially on $\Spin(\partial M_2)$:
this follows from the naturality of the Atiyah--Johnson correspondence 
$\Spin(\partial M_2) \to \Quad\left(\partial M_2\right)$ (see \S \ref{subsubsec:relative_Chern}).
According to Lemma \ref{lem:gluing}, there are two gluing maps
$$
\xymatrix{
\Spin^c(M_1,\sigma_1) \times \Spin^c(M_2,\sigma_2) 
\ar[r]^-{\cup_f} & \Spin^c(M)
}
$$
$$
\xymatrix{
\Spin^c(M_1,\sigma_1) \times \Spin^c(M_2,\sigma_2)
\ar[r]^-{\cup_{f\circ h}} & \Spin^c(M')
}
$$
which are affine, via Poincar\'e duality,
over $j_{1,*}\oplus j_{2,*}$ and $j'_{1,*}\oplus j'_{2,*}$ respectively. 
Since $\partial M_2$ is connected, the map $\cup_f$ is surjective.
\emph{Choose} $\alpha_1 \in \Spin^c(M_1,\sigma_1)$
and $\alpha_2 \in \Spin^c(M_2,\sigma_2)$ such that
$\alpha=\alpha_1\cup_f \alpha_2$, next set
\begin{displaymath} 
\alpha'=\alpha_1\cup_{f\circ h}\alpha_2 \in \Spin^c(M')
\end{displaymath}
and define $\Omega(\alpha)$ to be $\alpha'$.

We have to verify that $\Omega(\alpha)$
is well-defined by that procedure.
Assume other intermediate choices $\tilde{\sigma}_2$,
$\tilde{\alpha}_1$ and $\tilde{\alpha}_2$ instead of 
$\sigma_2$, $\alpha_1$ and $\alpha_2$ respectively, 
leading to $\tilde{\alpha}':=\tilde{\alpha}_1\cup_{f\circ h}\tilde{\alpha}_2$. 
We claim that $\alpha'=\tilde{\alpha}'$.

Consider first the particular case when 
$\tilde{\sigma}_2=\sigma_2 \in \Spin(\partial M_2)$. Since 
$\alpha_1\cup_f \alpha_2 = \alpha =
\tilde{\alpha}_1\cup_f \tilde{\alpha}_2$, we have that
\begin{displaymath}
j_{1,*} P^{-1}\left(\alpha_1-\tilde{\alpha}_1\right)
+ j_{2,*} P^{-1}\left(\alpha_2-\tilde{\alpha}_2\right)
= P^{-1}\left(\alpha-\alpha\right)=0\in H_1(M).
\end{displaymath}
Applying $\Phi$ to that identity, we obtain the equation
\begin{displaymath}
j'_{1,*} P^{-1}\left(\alpha_1 - \tilde{\alpha}_1\right)
+ j'_{2,*} P^{-1}\left(\alpha_2 - \tilde{\alpha}_2\right)
= 0\in H_1(M')
\end{displaymath}
whose left term equals $P^{-1}\left(\alpha'-\tilde{\alpha}'\right)$.
We conclude that $\alpha'=\tilde{\alpha}'$.

We now turn to the general case. For this, choose an arbitrary element
$$
\tau_2\in \Spinc\left([0,1]\times \partial M_2,0 \times (-\sigma_2) \cup 1 \times \tilde{\sigma}_2 \right).
$$
Having set $\tilde{\sigma}_1= f_*(-\tilde{\sigma}_2)$, define
$$
\tau_1=\left(\hbox{Id} \times f \right)_*(-\tau_2) 
\in \Spinc\left([0,1]\times \partial M_1,0 \times (-\sigma_1) \cup 1 \times \tilde{\sigma}_1 \right).
$$
Here, $-\tau_2\in \Spinc\left(-[0,1]\times \partial M_2,0 \times \sigma_2 \cup 1 \times (-\tilde{\sigma}_2) \right)$
is obtained from $\tau_2$ by time--reversing.
For $i=1,2$, the collar of $\partial M_i$ in $M_i$ and Lemma \ref{lem:gluing} give a map
$$
\xymatrix{\Spinc(M_i,\sigma_i) \times 
\Spinc\left([0,1]\times \partial M_i,0 \times (-\sigma_i) \cup 1 \times \tilde{\sigma}_i \right)
\ar[r]^-{\cup_{{\rm \footnotesize col}}} & \Spinc(M_i,\tilde{\sigma}_i)}.
$$
From the definition of the gluing map $\cup_f$ 
and by using the ``double collar'' of $\partial M_1\cong -\partial M_2$ in $M$, one sees that
$\alpha=\alpha_1\cup_f \alpha_2$ may also be written as
$$
\alpha=\left(\alpha_1 \cup_{{\rm \footnotesize col}} \tau_1 \right) \cup_f 
\left(\alpha_2 \cup_{{\rm \footnotesize col}} \tau_2\right).
$$
It follows from the special case treated previously that, whatever the choices
of $\tilde{\alpha}_1$ and $\tilde{\alpha}_2$ have been,
$$
\tilde{\alpha}'=\left(\alpha_1 \cup_{{\rm \footnotesize col}} \tau_1 \right) \cup_{f\circ h} 
\left(\alpha_2 \cup_{{\rm \footnotesize col}} \tau_2\right).
$$
On the other hand, having set 
$$
\tau_1'=\left(\hbox{Id} \times (f\circ h) \right)_*(-\tau_2) 
\in \Spinc\left([0,1]\times \partial M_1,0 \times (-\sigma_1) \cup 1 \times \tilde{\sigma}_1 \right),
$$
one sees that $\alpha'= \alpha_1 \cup_{f\circ h} \alpha_2$ may also be written as
$$
\alpha'=\left(\alpha_1 \cup_{{\rm \footnotesize col}} \tau'_1 \right) \cup_{f\circ h} 
\left(\alpha_2 \cup_{{\rm \footnotesize col}} \tau_2\right).
$$
Consequently, it is enough to prove that
\begin{equation}
\label{eq:basta}
\tau_1=\tau'_1 \in \Spinc\left([0,1]\times \partial M_1,0 \times (-\sigma_1) \cup 1 \times \tilde{\sigma}_1 \right).
\end{equation}
The latter space of relative \spinc--structures is classified by the Chern class map
since $H^2\left( [0,1]\times \partial M_1, \partial [0,1]\times \partial M_1\right)$ 
has no $2$--torsion. Moreover, the naturality of the Chern class 
and the fact that $h$ preserves the homology imply that
$$
c(\tau_1) = \left(\hbox{Id} \times f\right)_* (c(-\tau_2)) =
\left(\hbox{Id} \times (f\circ h) \right)_* (c(-\tau_2))= c(\tau_1').
$$
We conclude that identity (\ref{eq:basta}) holds and
that the map $\Omega$ is well-defined.

The fact that $\Omega$ is affine and the last statement
of the proposition are readily derived from the properties
of the gluing maps $\cup_f$ and $\cup_{f\circ h}$ stated in
Lemma \ref{lem:gluing}, and from the definition of the isomorphism $\Phi$.
\end{proof}

\begin{remark}
\label{rem:boundary_bis}
We could have considered as well the case
when $M_1$ and $M_2$ have disconnected boundary, 
but are glued together along a connected component of their boundary to give $M$
(so that $\partial M \cong \partial M' \neq \varnothing$).
Then, in view of Remark \ref{rem:boundary},
Proposition \ref{prop:twist} can easily be generalized
to involve \spinc--structures on $M$ and $M'$ relative to a fixed 
\spin--structure on their identified boundaries.
\end{remark}

\subsubsection{Definition of the $Y^c$--surgery move}

\label{subsubsec:def_Yc}
We  explain how $Y$--surgery makes sense in 
the setting of \spinc--manifolds. For simplicity, we consider only the case
of a closed oriented $3$--manifold $M$.

Let $j:H_3 \hookrightarrow M$ be an embedding. 
We denote by  $\Phi_j: H_1(M) \to H_1\left(M_j\right)$  the isomorphism
defined by the commutative diagram
$$
\xymatrix{
&H_1(M)\ar@{.>}[dd]^{\Phi_j}_\simeq\\
H_1\left(M\setminus {\rm{int}}\left({\rm{Im}}(j)\right) \right) \ar@{->>}[ru]^{k_*}
\ar@{->>}[rd]_{k'_*} & \\
&H_1\left(M_j\right)
}
$$
where $k:M\setminus {\rm{int}}\left({\rm{Im}}(j)\right) \hookrightarrow M$ and
$k':M\setminus {\rm{int}}\left({\rm{Im}}(j)\right) \hookrightarrow M_j$ denote inclusions.
\begin{lemma}
\label{lem:definition_Yc}
There exists a canonical bijection 
$$
\xymatrix{
{\rm{Spin}}^c(M) \ar[r]^-{\Omega_j}_-{\simeq} & {\rm{Spin}}^c\left(M_j\right), \quad
}
\alpha \longmapsto \alpha_j
$$ 
which is affine over $P\Phi_j P^{-1}$.
Moreover, the diagram 
\begin{displaymath}
\xymatrix{
{\rm{Spin}}^c(M)\ar[r]^-{\Omega_j}\ar[d]_c & {\rm{Spin}}^c(M_j)\ar[d]^c \\
H^2(M) \ar[r]_-{P \Phi_j P^{-1}} & H^2(M_j)
}
\end{displaymath}
is commutative.
\end{lemma}
\begin{proof}
By Remark \ref{rem:Y_as_a_twist}, one can define a
self-diffeomorphism $h$ of $\partial H_3$ 
acting trivially in homology and  such that there exists a diffeomorphism
\begin{displaymath}
\xymatrix{
{M_{j}= M \setminus \hbox{int}\left(\hbox{Im}(j)\right) 
\cup_{j|_{\partial H_3}} \left(H_{3}\right)_{B}} \ar[r]^-f_-\cong &
{M \setminus \hbox{int}\left(\hbox{Im}(j)\right)
\cup_{j|_{\partial H_3} \circ h} H_3}}
\end{displaymath}
which restricts to the identity on $M \setminus \hbox{int}\left(\hbox{Im}(j)\right)$.
This diffeomorphism induces a bijection
$$
\xymatrix{
{\Spin^c\left(M_j\right)} \ar[r]^-{f_*}_-\simeq &
{\Spin^c\left(M \setminus \hbox{int}\left(\hbox{Im}(j)\right) 
\cup_{j|_{\partial H_3} \circ h} H_3\right)}.
}
$$
Also, by \S \ref{subsubsec:twist}, there is a canonical bijection
$$
\xymatrix{
{\Spin^c\left(M\right)} \ar[r]^-{\Omega}_-\simeq &
{\Spin^c\left(M \setminus \hbox{int}\left(\hbox{Im}(j)\right)
\cup_{j|_{\partial H_3}\circ h} H_3\right)}.
} 
$$
We define $\Omega_j$ to be the composite $f_*^{-1}\Omega$.
This composite is easily verified to be independent of the pair $(h,f)$ with the above property.
\end{proof}

Let $G$ be a $Y$--graph in $M$.
Let also $j:H_3 \hookrightarrow M$ and $j':H_3 \hookrightarrow M$
be some trivializations of regular neighborhoods of $G$ in $M$.
There exists an ambiant isotopy $\left(q_t: M\to M\right)_{t\in [0,1]}$
between $j$ and $j'$: $q_0=\hbox{Id}_M$ and $q_1\circ j=j'$.
Let $q: M_j \to M_{j'}$ be the positive diffeomorphism induced by $q_1$
in the obvious way. One can verify that $q_*\circ \Omega_j = \Omega_{j'}$.
Thus, for any \spinc--structure $\alpha$ on $M$,
the \spinc--manifolds $\left(M_j,\alpha_j\right)$ and
$\left(M_{j'},\alpha_{j'}\right)$ are \spinc--diffeomorphic.
\begin{definition}
The Spin$^c$--manifold obtained from $(M,\alpha)$
by \emph{$Y^c$--surgery along $G$}, denoted by $\left(M_G,\alpha_G\right)$,
is the Spin$^c$--diffeomorphism
class of the manifold $\left(M_j,\alpha_j\right)$.
We call \emph{$Y^c$--equivalence} the equivalence relation 
among closed $3$--dimensional \spinc--manifolds
generated by $Y^c$--surgeries and \spinc--diffeomorphisms.
\end{definition}
\noindent
In the sequel, the notation $M_G$ will sometimes refer to
a representative $M_j$ obtained by fixing a trivialization $j$ of a
regular neighborhood of $G$ in $M$. 
Similarly, $\alpha_G$, $\Omega_G$ and $\Phi_G$ will stand for $\alpha_j$,
$\Omega_j$ and $\Phi_j$ respectively.
\begin{remark}
In the case of compact oriented $3$--manifolds with boundary,
the $Y^c$--surgery move is defined similarly using
\spinc--structures relative to \spin--structures.
(See  Remark \ref{rem:boundary_bis}.)
\end{remark}
It follows from the definition that, for any two disjoint $Y$--graphs $G_1$ and $G_2$  
in $M$, the \spinc--manifolds
$\left(\left(M_{G_1}\right)_{G_2},\left(\alpha_{G_1}\right)_{G_2}\right)$
and $\left(\left(M_{G_2}\right)_{G_1},\left(\alpha_{G_2}\right)_{G_1}\right)$ are \spinc--diffeomorphic.
So, the $Y^c$--surgery along a family of disjoint $Y$--graphs makes sense.
\begin{definition}
Let $I$ be an invariant of $3$--dimensional \spinc--manifolds
with values in an Abelian group $A$. The invariant $I$ is said
to be  \emph{of degree at most $d$} if, for any $3$--dimensional \spinc--manifold 
$(N,\sigma)$ and for any family $S$ of at least $d+1$ 
pairwise disjoint $Y$--graphs in $N$, the identity
\begin{equation}
\sum_{S'\subset S} (-1)^{|S'|}\cdot I(N_{S'},\sigma_{S'})=0 \in A
\end{equation}
holds. Here, the sum is taken over all sub-families $S'$ of $S$.
\end{definition}
\noindent
Thus, the $Y^c$--surgery move is the elementary move of a \spinc--refinement of 
the Goussa\-rov--Habiro theory of finite type invariants. In particular, two $3$--dimensional \spinc--man\-ifolds
are $Y^c$--equivalent if and only if they are not distinguished by degree $0$
invariants. It can be shown that the ``calculus of clovers'' from \cite{GGP}, 
which is equivalent to the ``calculus of claspers'' from \cite{Habiro}, extends to \spinc--manifolds.
\begin{remark} 
\label{rem:compatibility}
A \spin--refinement of the Goussarov--Habiro theory has been considered in \cite{Mas}. 
In particular, it is shown that the $Y$--surgery along $G$ induces a canonical bijection
$\Theta_G: {\Spin}(M) \to {\Spin}\left(M_G\right)$. Both refinements of the theory
are compatible, in the sense that the following diagram commutes:
$$
\xymatrix{
{\Spin}(M)\ar[d]_\beta \ar[r]^-{\Theta_G}_-{\simeq} 
& {\Spin}\left(M_G\right) \ar[d]^\beta\\
{\Spin}^c(M) \ar[r]^-{\simeq}_-{\Omega_G} & {\Spin}^c(M_G).
}
$$
\end{remark}

\subsubsection{A combinatorial description of the $Y^c$--equivalence relation}

A given equivalence relation among closed oriented $3$--manifolds
can sometimes be derived from an unknotting operation via surgery 
presentations in $\mathbf{S}^3$. It is well-known that the $Y$--equivalence relation 
can be formulated that way with the \emph{$\Delta$--move} of \cite{MN} as unknotting operation.
We refine this to the context of \spinc--manifolds.
\begin{lemma}
\label{lem:Yc_Deltac}
The $Y^c$--equivalence relation is generated
by \spinc--diffeomorphisms and $\Delta^c$--moves, 
if  the $\Delta^c$--move is defined to be the move
depicted on Figure \ref{fig:delta}
between surgery presentations of closed $3$--dimensional \spinc--manifolds (see \S\ref{subsubsec:combinatorial_spinc}).
\begin{figure}[h]
\begin{center}
\includegraphics[height=4cm,width=8cm]{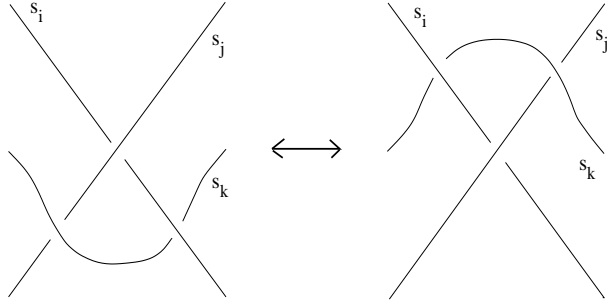}
\caption{A $\Delta^c$--move.}
\label{fig:delta}
\end{center}
\end{figure}
\end{lemma}
\begin{proof} 
Let $M$ be a closed connected oriented $3$--manifold 
and let $G$ be a $Y$--graph in $M$. Let $\psi: M\to V_L$  
be a surgery presentation of $M$, where $L$ 
is a $n$--component ordered oriented framed link in $\mathbf{S}^3$. 
Isotope $G$ in $M$ so that $\psi(G)$ becomes disjoint from the link dual to $L$, 
then $\psi(G)$ can be regarded as a subset of $\mathbf{S}^3\setminus L$.
In the image by $\psi$ of the regular neighborhood of $G$ in $M$,  put 
the $2$--component framed link $K$ depicted on Figure
\ref{fig:two-comp}. 
\begin{figure}[h]
\begin{center}
\includegraphics[height=4cm,width=5cm]{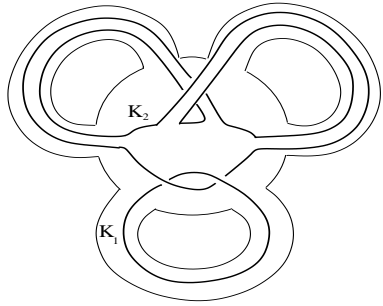}
\caption{$Y$--surgery as surgery along a $2$--component link.}
\label{fig:two-comp}
\end{center}
\end{figure}
The link $K$ can be obtained from the link $B$ of Figure \ref{fig:B}
by some slam dunks (see Example \ref{ex:slam_dunk}) 
and handle slidings in $H_3$. 
In particular, there is an obvious surgery presentation  
$\psi': M_G\to V_{L\cup K}$ induced by $\psi$.
With the viewpoint from \S \ref{subsubsec:combinatorial_spinc}, 
we want to identify the combinatorial analog of the bijection $\Omega_G$.
In other words, we look for the map $O_G$ making the diagram
$$
\xymatrix{
{\mathcal{V}_{L}} \ar@{.>}[r]^{O_{G}} & {\mathcal{V}_{L\cup K}} \\
{{\Spin}^c(V_L)}\ar[u]^{\simeq} &{{\Spin}^c(V_{L\cup K})}\ar[u]_{\simeq}\\
{{\Spin}^c(M)}\ar[u]^{\simeq}_{\psi_*} \ar[r]_{\Omega_G}& 
{{\Spin}^c\left(M_{G}\right)}\ar[u]_{\simeq}^{\psi'_*}
}
$$
commute. This is contained in the next claim, which will allow us to prove that 
the $\Delta^c$--move and the $Y^c$--surgery move are equivalent.
\begin{claim}
\label{claim:O_G}
Let $B_L$ denote the linking matrix of $L$ and let $K$ be 
appropriately oriented so that the ordered union
of ordered oriented framed links $L\cup K$ 
has its linking matrix of the form
\begin{displaymath}
B_{L\cup K}=\left(
\begin{array}{ccc|cc}
& & & x_{1} & 0\\
& B_{L} & & \vdots & \vdots\\
& & & x_{l} & 0\\
\hline
x_{1} &\cdots & x_{l} & x & 1\\
0 & \cdots & 0 & 1 & 0 
\end{array} \right).
\end{displaymath}
Then, the map $O_G$ sends a Chern vector $[s]$ 
to the Chern vector $\left[(s,x,0)\right]$.
\end{claim}
\begin{proof} 
As pointed out in Remark \ref{rem:compatibility},
a $Y$--surgery along $G$ induces a bijection
$\Theta_G:{\Spin}(M) \to {\Spin}(M_G)$,
a combinatorial analog of which is given in \cite{Mas}. 
Using the compatibility between $\Theta_G$ and $\Omega_G$ 
together with \S \ref{subsubsec:pont_combinatorial}, we
see that the claim holds at least for those Chern vectors
that come from $\mathcal{S}_L$.

Denote by $(H,f)$ the lattice corresponding to the intersection pairing
of $W_L$, and by $(H',f')$ that of $W_{L\cup K}$. 
Recall from Remark \ref{rem:combinatorial_Chern} that there are
canonical isomorphisms $H^2\left(V_L\right)\simeq 
\hbox{Coker}\ \widehat{f}$ and $H^{2}\left(V_{L\cup K}\right)
\simeq \hbox{Coker}\ \widehat{f'}$.
The isomorphism $P \Phi_G P^{-1}: H^2(M) \to H^2\left(M_G\right)$
corresponds then to the isomorphism $\hbox{Coker}\ \widehat{f}
\to \hbox{Coker}\ \widehat{f'}$ defined by $[y]\mapsto [\left(y,0,0\right)]$.

Take now $[s] \in \mathcal{V}_L$ arising
from $\mathcal{S}_L$ and let
$[y] \in\Z^n/\hbox{Im}\ B_L\simeq 
\hbox{Coker}\ \widehat{f}$. We aim to calculate 
$O_G\left([s]+ [y]\right) \in \mathcal{V}_{L\cup K}$.
The ``$+$'' here corresponds to the action of $H^2\left(V_L\right)$ on
${\Spin}^c(V_L)$ (see Remark \ref{rem:combinatorial_Chern}).
The map $\Omega_G$ being affine over $P \Phi_G P^{-1}$,
we have that $O_G\left([s] + [y]\right)$
$=O_G\left([s]\right) + \left[(y,0,0)\right]$
$= \left[(s,x,0)\right] + \left[(y,0,0)\right]$
$=\left[(s+2y,x,0)\right]$. Therefore, the claim
also holds for $[s] + [y]$$=[s+2y]$.
The transitivity of the action of $H^2(V_L)$ on ${\Spin}^c(V_L)$ allows
us to conclude.
\end{proof}

Figure \ref{fig:Y_to_delta} and Figure \ref{fig:delta_to_Y} prove that, 
up to \spinc--diffeomorphisms, a $\Delta^c$--move can be realized 
by a $Y^c$--surgery and vice versa.

In Figure \ref{fig:Y_to_delta}, the first \spinc--diffeomorphism is obtained 
by applying Claim \ref{claim:O_G}, while the second one is obtained from
one handle sliding and one slam dunk.
\begin{figure}[h!]
\begin{center}
\includegraphics[width=11cm,height=4cm]{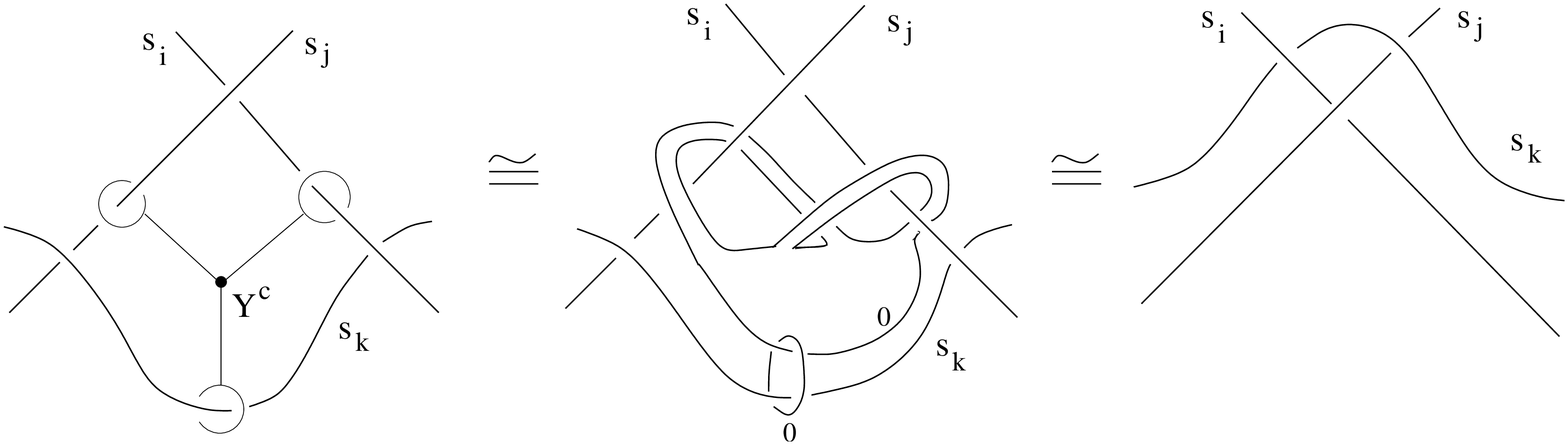}
\caption{A $\Delta^c$--move can be realized by a $Y^c$--surgery.}
\label{fig:Y_to_delta}
\end{center}
\end{figure}

In Figure \ref{fig:delta_to_Y}, the first \spinc--diffeomorphism is obtained
from three slam dunks.
Next, a $\Delta^c$--move is applied. The second \spinc--diffeomorphism
is obtained by \spinc{} Kirby's calculi (in particular, two slam dunks have been performed),
and the last one is obtained from Claim \ref{claim:O_G}.
\begin{figure}[h!]
\begin{center}
\includegraphics[width=10cm,height=6cm]{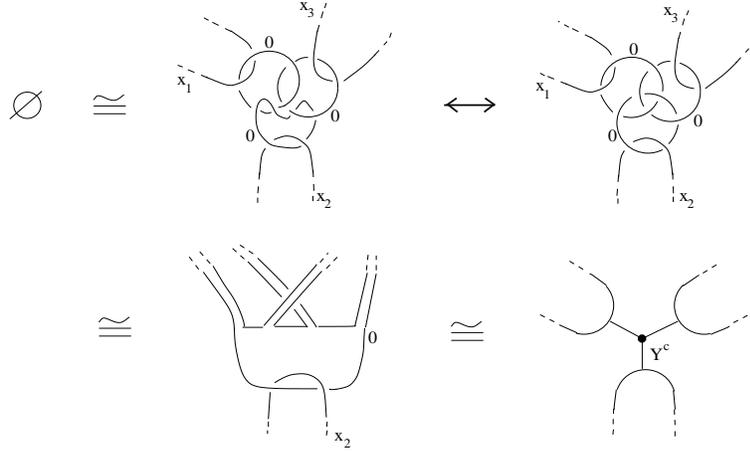}
\caption{A $Y^c$--surgery can be realized by a $\Delta^c$--move.}
\label{fig:delta_to_Y}
\end{center}
\end{figure}
\end{proof}

\subsection{Proof of Theorem \ref{th:Matveev_spinc}}

In this subsection, we prove the characterization of the $Y^c$--equivalence relation, as
announced in the introduction. We need two results  
concerning classification of quadratic functions up to isomorphism, proved in \cite{DM1}.

\subsubsection{Isomorphism classes of quadratic functions}

\label{subsubsec:iso_classes}
There is a natural notion of isomorphism among triples $(H,f,c)$ 
defined by bilinear lattices with characteristic form 
(see \S\ref{subsec:discriminant}): we say that two triples 
$(H,f,c)$ and $(H',f',c')$ are \emph{isomorphic} if there is an isomorphism 
$\psi:H \to H'$ such that $f = f' \circ (\psi \times \psi)$ 
and $c= c' \circ \psi$ mod $2\widehat{f}(H)$.  
Such triples form a monoid for the orthogonal sum 
$\oplus$. Two triples $(H,f,c)$ and $(H',f',c')$ are said to be 
\emph{stably equivalent} if they become isomorphic after stabilizations 
with some copies of $(\Z,\pm 1,\hbox{Id})$, which denotes 
the bilinear lattice defined on $\Z$ by $(1,1)\mapsto \pm 1$ and equipped 
with the characteristic form $\hbox{Id}=\hbox{Id}_{\Zi}$.\\
Note that, for any bilinear lattices $(H,f)$ and $(H',f')$, there is a map 
$$
\psi \mapsto \psi^{\sharp}, \ \ \hbox{Iso}\left(\hbox{Coker}\ \widehat{f},  \hbox{Coker}\ \widehat{f'}\right)
\to  \hbox{Iso}\left(G_{f'}, G_{f}\right)
$$
since the pairing (\ref{eq:evaluation}) is right nonsingular.
\begin{theorem} \cite{DM1}
\label{th:stabilisation}
Two bilinear lattices with characteristic form $(H,f,c)$ and $(H',f',c')$ 
are stably equivalent if, and only if, there exists  an element 
$$\psi^{\sharp} \in
{\rm{Im}} \left({\rm{Iso}}\left({\rm{Coker}}\ \widehat{f},
{\rm{Coker}}\ \widehat{f'}\right) \to
{\rm{Iso}}\left(G_{f'}, G_{f}\right) \right)$$ 
such that the associated
quadratic functions $(G_f,\phi_{f,c})$ and $(G_{f'},\phi_{f',c'})$
are isomorphic via $\psi^{\sharp}$.
Furthermore, any such isomorphism  between $(G_{f'},\phi_{f',c'})$ and
$(G_{f},\phi_{f,c})$ lifts to a stable equivalence 
between $(H,f,c)$ and $(H',f',c')$.
\end{theorem}
\begin{remark}
\label{rem:reconnaissance_image}
Let $\Psi$ be an  isomorphism  between $(G_{f'},\phi_{f',c'})$ and
$(G_{f},\phi_{f,c})$ and suppose that $f$ and $f'$ are degenerate. 
Then, $\Psi$ does not necessarily arise from an isomorphism
$\psi:{\rm{Coker}}\ \widehat{f} \to {\rm{Coker}}\ \widehat{f'}$.
In fact, it does if and only if $\Psi|_{{\rm{\scriptsize Ker}}\ \widehat{L_{f'}}}: 
\hbox{Ker}\ \widehat{L_{f'}} \to \hbox{Ker}\ \widehat{L_{f}} $ lifts to an isomorphism
${\rm{Ker}}\ \widehat{f'} \to {\rm{Ker}}\ \widehat{f}$. (See \cite{DM1} for details.)
\end{remark}

Let now $q: G \to \Q/\Z$ be a quadratic function on an Abelian group $G$. 
We shall say that $q$ meets the {\emph{finiteness condition}} if 
\begin{enumerate}
\item[$\centerdot$] $G/\hbox{Ker}\ \widehat{b_q}$ is finite;
\item[$\centerdot$] the extension $G$ of $\hbox{Ker}\ \widehat{b_q}$ by 
$G/\hbox{Ker}\ \widehat{b_q}$ is split.
\end{enumerate}
We shall also denote by $r_q$ the homomorphism obtained by restricting
$q$ to $\hbox{Ker}\ \widehat{b_q}$. 
\begin{theorem}\cite{DM1}
\label{th:invcomplet}
Two quadratic functions $q:G\to \Q/\Z$ and $q':G'\to \Q/\Z$ 
satisfying the finiteness condition are isomorphic 
if, and only if, there is an isomorphism $\Psi:G'\to G$ such that 
$b_{q'} = b_{q}\circ(\Psi \times \Psi)$, $d_{q'} = d_{q}\circ \Psi$, 
$r_{q'} = r_{q}\circ\Psi|$ and $\gamma(q' \circ s') = \gamma(q \circ s)$ 
for some $\Psi$--compatible sections $s$ and $s'$ of the canonical
epimorphisms $G\to  G/{\rm{Ker}}\ \widehat{b_q}$
and $G'\to  G'/{\rm{Ker}}\ \widehat{b_{q'}}$.
\end{theorem}
\noindent
Here, the $\Psi$--compatibility condition refers to the commutativity of the diagram
\begin{displaymath}
\xymatrix{
{G'} \ar[d]_-{\Psi}^-{\simeq} &
{G'/{\rm{Ker}}\ \widehat{b_{q'}}} \ar[l]_-{s'} \ar[d]^-{[\Psi]}_-{\simeq} \\
{G} & {G/{\rm{Ker}}\ \widehat{b_{q}}} \ar[l]^-{s}
}
\end{displaymath}
where $[\Psi]$ is the isomorphism induced by $\Psi$.
\begin{remark} 
\label{rem:subtilite}
Theorem \ref{th:invcomplet} does not claim that 
$q'=q\circ \Psi$ if the four conditions hold. 
Nevertheless, as follows from the proof in \cite{DM1}, it is true
that there exists an isomorphism $\varphi:G' \to G$ such that 
$q'=q\circ \varphi$ and $\varphi|_{{\rm{Ker}}\ \widehat{b_{q'}}}=
\Psi|_{{\rm{Ker}}\ \widehat{b_{q'}}}$.
\end{remark}
We now go into the proof of Theorem \ref{th:Matveev_spinc}.
In the sequel, we consider
two closed connected $3$--dimensional  \spinc--manifolds,
$(M,\sigma)$ and $(M',\sigma')$.

\subsubsection{Proof of the equivalence $(2) \Longleftrightarrow (3)$
of Theorem \ref{th:Matveev_spinc}}

Next lemma is easily proved from the definitions.
\begin{lemma} 
\label{lem:equivalences}
Let $\psi: H_1(M) \to H_1(M')$ be an isomorphism, which induces
a dual isomorphism $\psi^\sharp:$ $H_2(M';\Q/\Z) \to H_2(M;\Q/\Z)$ 
with respect to the intersection pairings. The following assertions are equivalent:
\begin{enumerate}
\item[(a)] $L_{M'}=L_{M}\circ \left(\psi^{\sharp}\times \psi^\sharp\right)$.
\item[(b)] $\lambda_{M}=\lambda_{M'} \circ \left(\psi|\times \psi|\right)$.
\item[(c)] The following diagram is commutative:
$$
\xymatrix{
H_2\left(M';\Q/\Z\right) \ar[d]_-{\psi^{\sharp}}^{\simeq}
\ar[r]^-{B} & {\rm{Tors}}\ H_1(M')\\
H_2\left(M;\Q/\Z\right) \ar[r]^-{B} & {\rm{Tors}}\ H_1(M). \ar[u]^-{\psi|}_{\simeq}
}
$$
\end{enumerate}
\end{lemma}

Suppose that the condition (2) of Theorem \ref{th:Matveev_spinc} is satisfied.
This implies  that $L_{M'}=L_{M}\circ \left(\psi^{\sharp}\times\psi^\sharp\right)$ 
and so that $\lambda_{M}=\lambda_{M'} \circ \left(\psi|\times \psi|\right)$
by Lemma \ref{lem:equivalences}.

Condition (2) also implies the relation $d_{\phi_{M',\sigma'}}=
d_{\phi_{M,\sigma}}\circ \psi^\sharp$ between homogeneity defects of 
quadratic functions. So, by Lemma \ref{lem:homogeneity-defect}, we have 
$\langle c(\sigma'),x'\rangle= \langle c(\sigma),\psi^\sharp(x')\rangle$
for all $x'\in H_2(M';\Q/\Z)$. By left nondegeneracy of the pairing
$\bullet: H_1(M')\times H_2(M';\Q/\Z) \to\Q/\Z$, we conclude that
$P^{-1}c(\sigma')=\psi\left( P^{-1}c(\sigma)\right)$.

Last, the quadratic function
$$\phi_{M,\sigma} \circ s=\phi_{M,\sigma} \circ \psi^\sharp \circ s' \circ \psi|=
\phi_{M',\sigma'} \circ s'\circ \psi|$$
is isomorphic to $\phi_{M',\sigma'}\circ s'$: hence, these two quadratic
functions have identical Gauss sums. Therefore the condition (3) holds.\\

Conversely, suppose that the condition (3) of Theorem \ref{th:Matveev_spinc} 
is satisfied. The short exact sequence
\begin{displaymath}
\xymatrix{
{0}\ar[r]& {H_2(M)\otimes \Q/\Z}\ar[r]& {H_2(M;\Q/\Z)} \ar[r]^B
& {\hbox{Tors}\ H_1(M)} \ar[r] & {0} 
}
\end{displaymath}
is split, we have that $H_2(M)\otimes \Q/\Z=\hbox{Ker}\ \widehat{L_M}$ 
and $\hbox{Tors } H_1(M)$ is finite: 
thus, $\phi_{M,\sigma}$ 
meets the finiteness condition of \S \ref{subsubsec:iso_classes}.
Since $\lambda_{M}= \lambda_{M'}\circ \left(\psi|\times \psi|\right)$, we obtain 
by Lemma \ref{lem:equivalences} that 
$L_{M'}=L_{M}\circ \left(\psi^{\sharp}\times \psi^\sharp\right)$. 
Since $\psi\left(P^{-1}c(\sigma)\right) = P^{-1}c\left(\sigma'\right)$, 
we deduce from Lemma \ref{lem:behaviour_on_the_kernel} and
Lemma \ref{lem:homogeneity-defect} that
$r_{\phi_{M',\sigma'}} = r_{\phi_{M,\sigma}} \circ \psi^\sharp|$ and that
$d_{\phi_{M',\sigma'}} = d_{\phi_{M,\sigma}} \circ \psi^\sharp$ respectively.
Also, since $\psi| \circ B \circ \psi^\sharp=B$ (by Lemma \ref{lem:equivalences}), 
the $\psi$--compatibility condition between $s$ and $s'$ required by the condition (3) of
Theorem \ref{th:Matveev_spinc} coincides with the $\psi^\sharp$--compatibility in the sense of 
\S \ref{subsubsec:iso_classes}. Therefore, by Theorem \ref{th:invcomplet}, 
the quadratic functions $\phi_{M,\sigma}$ and $\phi_{M',\sigma'}$ are isomorphic.
More precisely, according to Remark \ref{rem:subtilite}, there exists an isomorphism
$\varphi: H_2(M';\Q/\Z) \to H_2(M;\Q/\Z)$ such that 
$\phi_{M',\sigma'}=\phi_{M,\sigma}\circ \varphi$ and $\varphi|_{H_2(M')\otimes \Qi/\Zi}$
coincides with $\psi^\sharp|_{H_2(M')\otimes \Qi/\Zi}=
\psi^{\sharp}\otimes \Q/\Z$. This latter fact, together with 
Remark \ref{rem:reconnaissance_image}, allows us to precise that
$\varphi$ equals $\eta^\sharp$ for a certain isomorphism 
$\eta: H_1(M) \to H_1(M')$. Consequently,
$\phi_{M',\sigma'}=\phi_{M,\sigma}\circ \eta^\sharp$.

\subsubsection{Proof of the equivalence $(1) \Longleftrightarrow (2)$ of Theorem \ref{th:Matveev_spinc}}

We prove implication $(1) \Longrightarrow (2)$ first. By Lemma
\ref{lem:Yc_Deltac}, it suffices to prove it when 
$(M,\sigma)$ and $(M',\sigma')$ are related by one \spinc--diffeo\-morphism
or, for some fixed surgery presentations, by one $\Delta^c$--move.
The first case follows immediately from the definition of the linking quadratic function.
The second case is deduced from the combinatorial formula
for the latter given at the end of \S \ref{subsec:4def}, 
and from the fact that a $\Delta$--move between ordered oriented framed links preserve the linking matrices.\\

Suppose now that condition (2) is satisfied. We can assume that $M=V_L$ and $M'=V_{L'}$,
where $L$ and $L'$ are ordered oriented framed links in $\mathbf{S}^3$. 
As in \S \ref{subsec:combinatorial}, we denote by $(H,f)$ and $(H',f')$
the intersection pairings of $W_L$ and $W_{L'}$ respectively. Let also $c\in \hbox{Char}(f)$ and
$c'\in \hbox{Char}(f')$ represent $\sigma$ and $\sigma'$ respectively. 
By hypothesis, the quadratic functions $\phi_{f,c}:G_f\to \Q/\Z$ and 
$\phi_{f',c'}:G_{f'}\to \Q/\Z$ are isomorphic via an isomorphism 
which is induced by an isomorphism $\hbox{Coker}\ \widehat{f}\to 
\hbox{Coker}\ \widehat{f'}$. So, by Theorem \ref{th:stabilisation}, the
bilinear lattices with characteristic form $(H,f,c)$ and $(H',f',c')$ are stably equivalent.

An isomorphism of bilinear lattices with characteristic form can be topologically realized 
by a finite sequence of Spin$^c$ Kirby's moves (see Theorem \ref{th:Spin^c_Kirby}): 
handle slidings and reversings of orientation.
Similarly, a stabilization by $(\Z,\pm 1,\hbox{Id})$
corresponds to a stabilization by the unknot.
Therefore, we can suppose, without loss of generality, that 
$(H,f,c) \simeq (H',f',c')$ through the isomorphism that identifies the preferred
basis of $H$ with that of $H'$. Concretely, this means that the linking matrices 
$B_L$ and $B_{L'}$ are equal and that there is a multi-integer $s$
such that the Chern vectors $[s]\in \mathcal{V}_L$ and $[s]\in \mathcal{V}_{L'}$
represent $\sigma$ and $\sigma'$ respectively.

A theorem\footnote{In fact, the first reference is \cite{Matveev} but the proof there is not detailed.}
of Murakami and Nakanishi \cite[Theorem 1.1]{MN} states that two ordered oriented framed  links 
have identical linking matrices if, and only if, they are $\Delta$--equivalent. 
Then, the ``decorated'' links $(L,s)$ and $(L',s)$ are $\Delta$--equivalent: 
therefore, by Lemma \ref{lem:Yc_Deltac}, the \spinc--manifolds $(M,\sigma)$ and $(M',\sigma')$ are $Y^c$--equivalent.
\begin{remark}
Observe that the present proof allows for a more precise statement 
of the equivalence $(1) \Longleftrightarrow (2)$ 
of Theorem \ref{th:Matveev_spinc}. Any finite sequence of \spinc--diffeomorphisms
and $Y^c$--surgeries
$$(M,\sigma)=(M_0,\sigma_0)\leadsto(M_1,\sigma_1)\leadsto (M_2,\sigma_2) 
\leadsto \cdots\leadsto (M_n,\sigma_n)=(M',\sigma')$$
yields an isomorphism $\psi: H_1(M) \to H_1(M')$. This is the composite
of the isomorphisms $H_1(M_i)\to H_1(M_{i+1})$, which is taken to be either
$g_*$ if the step $(M_i,\sigma_i)\leadsto (M_{i+1},\sigma_{i+1})$ is a 
\spinc--diffeomorphism $g$, either the isomorphism $\Phi_G$ 
if the step is the $Y^c$--surgery along a $Y$--graph $G\subset M_i$ (\S \ref{subsubsec:def_Yc}).
This isomorphism $\psi$ satisfies $\phi_{M',\sigma'}=\phi_{M,\sigma}\circ \psi^\sharp$. 
Conversely, given an isomorphism $\psi: H_1(M) \to H_1(M')$ with this property, 
one can find a finite sequence of \spinc--diffeomorphisms and $Y^c$--surgeries
from $(M,\sigma)$ to $(M',\sigma')$ inducing $\psi$ at the level of $H_1(-)$.
Here, we use the second statement of Theorem \ref{th:stabilisation}.
\end{remark}

\subsection{Applications and problems}

We conclude this paper with some applications of our results
illustrated by a few examples. We also state a few problems.

\subsubsection{The quotient set ${\rm{Spin}}^c(M)/Y^c$} 

Given a closed oriented $3$--manifold $M$, one may consider
the quotient set $$\hbox{Spin}^c(M)/Y^c$$
of \spinc--structures on $M$ modulo the $Y^c$--equivalence relation.
Let us consider a few examples. 
\begin{example} Take $M=\R \hbox{P}^3$. This manifold has 
two distinct \spinc--structures $\sigma_0$ and $\sigma_1$,
both arising from  \spin--structures. The quadratic functions 
$\phi_{M,\sigma_0}$ and $\phi_{M,\sigma_1}$ have different Gauss
sums (which are $\exp(2i\pi/8)$ and $\exp(-2i\pi/8)\in \C$). 
Therefore, by Corollary \ref{cor:QHS}, $\sigma_0$ is not
$Y^c$--equivalent to $\sigma_1$.
\end{example}
\begin{example} Take $M$ such that $H_1(M)\simeq \Z^n$. According 
to Corollary \ref{cor:torsion_free}, the set ${\hbox{Spin}^c(M)/Y^c}$ can be 
identified with $\left(2\Z^n\right) / \hbox{GL}(n;\Z)$ by the Chern class map.

In particular, if $M=\mathbf{S}^2\times \mathbf{S}^1$ and if an isomorphism
$H_1\left(M\right)\simeq \Z$ is fixed,
we denote by $\alpha_k$ the unique element of $\hbox{Spin}^c(M)$
such that $c(\alpha_k)=2k\in \Z$, with $k\in \Z$. Then, the $Y^c$--equivalence
classes are $\{\alpha_0\}$ and $\{\alpha_k,\alpha_{-k}\}$ with $k>0$. 
Observe from Theorem \ref{th:Spin^c_Kirby}, that these classes 
coincide with the diffeomorphism classes.
\end{example}

\begin{example} Take $M=\left(\mathbf{S}^2\times \mathbf{S}^1\right)\sharp\ \R \hbox{P}^3$.
By applying equivalence $(1) \Leftrightarrow (2)$ of Theorem \ref{th:Matveev_spinc}, 
the $Y^c$--equivalence classes are seen to be 
$\{\alpha_0\sharp \sigma_0\}$,  $\{\alpha_0\sharp \sigma_1\}$,
$\{\alpha_k\sharp \sigma_0,$ $ \alpha_k\sharp \sigma_1,$ $ \alpha_{-k}\sharp \sigma_0,$ 
$\alpha_{-k}\sharp \sigma_1 \}$ with $k>0$ odd, 
$\{\alpha_k\sharp\sigma_0,\alpha_{-k}\sharp \sigma_0\}$ 
and $\{\alpha_k\sharp\sigma_1,\alpha_{-k}\sharp \sigma_1\}$ with $k>0$ even.
Again, observe from Theorem \ref{th:Spin^c_Kirby}, that these classes 
coincide with the diffeomorphism classes.
\end{example}
In light of the previous examples, it is natural to ask whether the diffeomorphism classes of 
\spinc--structures of a given closed oriented $3$--manifold $M$ coincide with the
$Y^c$--equivalence classes. To answer this question by the negative, 
let us consider a class of manifolds for which
the \spinc--structures have been classified: the family of lens spaces.
Let $p\geq 2$ be an integer, let $q_1$, $q_2$ be some invertible elements of $\Z_p$
and let $L(p;q_1,q_2)$ be the corresponding lens space with the orientation induced 
from the canonical orientation of $\mathbf{S}^3$.
\begin{theorem} {\rm{(Turaev \cite{TReid})}}
\label{th:classification}
The number of orbits of \spinc--structures under the action of 
the group of positive self-diffeomorphisms of $L(p;q_1,q_2)$ is 
\begin{itemize}
\item[$\centerdot$] $\left[p/2\right] +1$, if $q_1^2\neq q_2^2$ or $q_1=\pm q_2$,
\item[$\centerdot$] $p/2 - b(p;q_1,q_2)/4+c(p;q_1,q_2)/2$, if $q_1^2=q_2^2$ and $q_1\neq \pm q_2$.
\end{itemize}
Here, for $x\in \Q$, $[x]$ denotes the greatest integer less or equal than $x$, $b(p;q_1,q_2)$ 
is the number of $i\in \Z_p$ for which $i$, $q_1+q_2-i$ and $q_2q_1^{-1}i$ are pairwise different, and 
$c(p;q_1,q_2)$ is the number of $i\in \Z_p$ such that $i=q_1+q_2-i=q_2q_1^{-1}i$.
\end{theorem}
\begin{proof} 
In \cite[\S 9.2.1]{TReid}, the Euler structures on $L(p;q_1,q_2)$ are classified up
to diffeomorphisms. The same kind of arguments can be used to classify these
up to \emph{positive} diffeomorphisms. Details are left to the reader.
\end{proof}

\noindent
The classification of the \spinc--structures on $L(p;q_1,q_2)$ up to 
$Y^c$--equivalence is easily obtained from Corollary \ref{cor:QHS}. For instance,
let us suppose that $p$ is odd. Then, $\hbox{Spin}^c\left(L(p;q_1,q_2)\right)/Y^c$ 
can be identified via the Chern class map with the quotient set $\Z_p/\sim$, where
$$
\forall i,j \in \Z_p,\quad \left(i\sim j\right) \iff 
\left(\exists r \in \Z_p,\ r^2=1\ \hbox{and}\ j=r i \right).
$$ 

\begin{example} Let $k\geq 4$ be an even integer and let $p=k^2-1$. 
Then, there are some \spinc--structures on $L(p;1,1)$ which are $Y^c$--equivalent
but which are not diffeomorphic. Indeed, according to Theorem \ref{th:classification}, 
$\hbox{Spin}^c\left(L(p;1,1)\right)$ contains $\left(p-1\right)/2+1$ diffeomorphism classes.
But, $k^2= 1\in \Z_p$ and $k\neq \pm 1 \in \Z_p$, so the cardinality of 
$\hbox{Spin}^c(L(p;1,1))/Y^c$ is strictly less than $(p-1)/2+1$.
\end{example}

\subsubsection{Reidemeister--Turaev torsions}

Let $\tau(M,\sigma)$ denote the maximal Abelian Reid\-emeister--Turaev torsion 
of a closed oriented $3$--manifold $M$ equipped with an Euler structure or,
equivalently, a \spinc--structure $\sigma$ \cite{TBook}. 
If $M$ is a rational homology sphere,
it turns out that $\phi_{M,\sigma}$ can be explicitely computed from $\tau(M,\sigma)$ \cite{Nic,DM2}. 
Thus, according to Corollary \ref{cor:QHS}, part of $\tau(M,\sigma)$ is of degree $0$.

\begin{problem}
\label{quest:torsions}
Derive from Reidemeister--Turaev torsions higher degree finite
type invariants of closed $3$--dimensional \spinc--manifolds.
\end{problem}

In the last chapter of \cite{MasThese}, it is studied how  Reidemeister--Turaev torsions 
vary under those twists defined in \S \ref{subsubsec:twist}. 
This variation is difficult to control for a generic $Y$--graph. Nevertheless, this variation
can be calculated explicitely in case of ``looped clovers''.
It is shown that Reidemeister--Turaev torsions satisfy a certain multiplicative degree $1$ relation
involving surgeries along looped clovers.

\subsubsection{From the \spin--refinement of the theory to its \spinc--refinement}

According to Remark \ref{rem:compatibility}, any \spinc--invariant of degree $d$ 
in the Goussarov--Habiro theory induces a \spin--invariant of degree $d$. 
The converse is not true.

\begin{example}
The Rochlin invariant $R(M,\sigma)\in \Z_{16}$ of a closed \spin--manifold $(M,\sigma)$ of
dimension $3$ is a finite type invariant of degree $1$ \cite{Mas}. But, it does not lift 
to an invariant of \spinc--manifolds in general. Indeed, consider the torus $\mathbf{T}^3$ and 
its canonical \spin--structure $\sigma^0$ (induced by its Lie group structure),
choose also $\sigma'$ in $\Spin\left(\mathbf{T}^3\right)$ different from $\sigma^0$. Then,
$\beta(\sigma')$ and $\beta\left(\sigma^0\right)$ coincide,
but $R\left(\mathbf{T}^3,\sigma^0\right)=8$ is not equal to 
$R\left(\mathbf{T}^3,\sigma'\right)=0$.
\end{example}

On the contrary, we have in degree $0$ the following consequence of both Theorem 
\ref{th:Matveev_spinc} and \cite[Theorem 1]{Mas}.
\begin{corollary}
Let $(M,\sigma)$ and $(M',\sigma')$
be closed $3$--dimensional \spin--manifolds. 
Then, $(M,\sigma)$ and $(M',\sigma')$ are 
distinguished by degree $0$ \spin--invariants 
if and only if $\left(M,\beta(\sigma)\right)$ and $\left(M',\beta(\sigma')\right)$
are distinguished by degree $0$ \spinc--invariants.
\end{corollary}

\begin{problem}
Compare in higher degrees the \spinc--refinement of the Goussarov--Habiro theory
with its \spin--refinement. 
\end{problem}

\bibliographystyle{amsalpha}

\vspace{1cm}

\textsc{\footnotesize F. Deloup, Einstein Institute of Mathematics,
Edmond J. Safra Campus -- Givat Ram,
The Hebrew University of Jerusalem,
91904 Jerusalem, Israel.}\\
\emph{\footnotesize and} \textsc{\footnotesize Laboratoire Emile Picard
(UMR 5580 CNRS/Universit\'e Paul Sabatier),
118 route de Narbonne, 31062 Toulouse Cedex 04, France.}

\emph{\footnotesize E-mail addresses:}
\texttt{\footnotesize deloup@math.huji.ac.il} \emph{\footnotesize \ and \ }
\texttt{\footnotesize deloup@picard.ups-tlse.fr}

\vspace{0.5cm}

\textsc{ \footnotesize G. Massuyeau, Laboratoire Jean Leray
(UMR 6629 CNRS/Universit\'e de Nantes),
2 rue de la Houssi\-ni\`e\-re, BP 92208,  44322 Nantes Cedex 03, France.}

\emph{\footnotesize E-mail address:} \texttt{\footnotesize 
massuyea@math.univ-nantes.fr}\\

\end{document}